\documentclass[oneside,reqno,11pt,a4paper]{amsart}
\usepackage[T1]{fontenc}
\usepackage[margin=2.5cm]{geometry}
\usepackage{amsmath,amsfonts,amssymb,amscd,amsthm}
\usepackage[usenames,dvipsnames,svgnames]{xcolor}
\usepackage[utf8x]{inputenc}
\usepackage{graphicx}
\graphicspath{{}{graphs/}{}}
\usepackage{caption}
\usepackage{subcaption}
\usepackage[numbers]{natbib}
\usepackage{doi}
\usepackage{autonum}
\usepackage{makecell}
\usepackage{hyperref}
\usepackage{acronym}
\usepackage{algorithm,algorithmic}
\usepackage{textgreek}
\usepackage{url}
\usepackage{color}
\usepackage{framed}
\usepackage{verbatim}
\usepackage{fancyvrb}
\usepackage{calc}
\usepackage{newtxtext}
\usepackage{newtxmath}
\usepackage{microtype} 
\usepackage[final]{pdfpages}

\usepackage{ragged2e}

\usepackage{booktabs}
\usepackage{tikz}
\definecolor{myellow}{RGB}{255,230,128}
\definecolor{gray20}{RGB}{204,204,204}
\definecolor{mygray}{RGB}{204,204,204}
\definecolor{mygreen}{RGB}{138,203,95}
\definecolor{myblue}{RGB}{77,151,214}

\newif\iflistings
\listingsfalse

\iflistings
   \usepackage{minted}

\newcommand{\julia}[1]{\mint{julia}|#1|}
\fi

\hypersetup{
    breaklinks=true,
    bookmarksopen=true,
    pdftitle={Intersecting surface meshes and Cartesian background meshes for embedded simulations},    
    pdfauthor={Santiago Badia, P. A. Martorell, F. Verdugo},     
    colorlinks=true,       
    linkcolor=black,          
    citecolor=blue,        
    filecolor=black,      
    urlcolor=blue           
}

\definecolor{bg}{rgb}{0.93,0.93,0.93}

\newtheorem{theorem}{Theorem}[section]

\theoremstyle{definition}
\newtheorem{definition}[theorem]{Definition}

\acrodef{amr}[AMR]{adaptive mesh refinement}
\acrodef{cad}[CAD]{computer-aided design}
\acrodef{cae}[CAE]{computer-aided engineering}
\acrodef{pde}[PDE]{partial differential equation}
\acrodefplural{pde}[PDEs]{partial differential equations}
\acrodef{cae}[CAE]{computer-aided engineering}
\acrodef{stl}[STL]{Stereolithography}
\acrodef{fe}[FE]{finite element}
\acrodef{vof}[VOF]{volume of fluid}
\acrodef{fem}[FEM]{finite element method}
\acrodef{xfem}[XFEM]{extended finite element method}
\acrodef{dof}[DOF]{degree of freedom}
\acrodefplural{dof}[DOFs]{degrees of freedom}
\acrodef{agfe}[agFE]{aggregated finite element}
\acrodef{agfem}[AgFEM]{aggregated finite element method}
\acrodef{cg}[CG]{continuous Galerkin}
\acrodef{dg}[DG]{discontinuous Galerkin}
\acrodef{amg}[AMG]{Algebraic Multi-Grid}
\acrodef{lic}[LIC]{line integral convolution}

\newcommand{\fig}[1]{Figure~\ref{#1}}
\newcommand{\sect}[1]{Section~\ref{#1}}
\newcommand{\alg}[1]{Algorithm~\ref{#1}}
\newcommand{\tab}[1]{Table~\ref{#1}}

\def\C{\mathcal{C}}
\def\esn{\epsilon_{\mathrm{sn}}}
\def\ehs{\epsilon_{\mathrm{hs}}}
\def\T{\mathcal{T}}
\def\B{\mathcal{B}}
\def\Tcut{\mathcal{T}^{\mathrm{cut}}}
\def\Bcut{\B^{\mathrm{cut}}}

\definecolor{shadecolor}{gray}{.92}
\definecolor{incolor}{rgb}{0,0,.7}
\definecolor{outcolor}{rgb}{.65,0,0}
\definecolor{syntaxcolor}{rgb}{.65,0,0}

\newif\ifsvgs
\svgsfalse

\newcommand{\includefig}[3][\tiny]{%
    \def\svgwidth{#2}
    #1
    
  \ifsvgs
    \updatepdffromsvg{#3}
  \fi
  \input{#3.pdf_tex}

}

\newcommand{\updatepdffromsvg}[1]{
  \executeiffilenewer{#1.svg}{#1.pdf}%
  {inkscape -z -C --file=#1.svg %
  --export-pdf=#1.pdf --export-latex && %
  sed "s|$(basename #1).pdf|$(echo #1).pdf|g" -i #1.pdf_tex }%
}

\newcommand{\executeiffilenewer}[3]{%
  \IfFileExists{#2}{}{\immediate\write18{#3}}
  \ifnum\pdfstrcmp{\pdffilemoddate{#1}}%
  {\pdffilemoddate{#2}}>0%
  {\immediate\write18{#3}}\fi%
}

\begin{document}

\title[Geometrical discretisation for unfitted finite elements on explicit representations]{Geometrical discretisations for unfitted finite elements on explicit boundary representations}

\author[S. Badia]{Santiago Badia$^{1,3,*}$}

\author[P. Martorell]{Pere A. Martorell$^{2}$}

\author[F. Verdugo]{Francesc Verdugo$^{3}$}

\thanks{\null\\
$^{1}$ School of Mathematics, Monash University, Clayton, Victoria, 3800, Australia.\\
$^{2}$ Department of Civil and Environmental Engineering, Universitat Polit\`ecnica de Catalunya, Edifici C2, Campus Nord UPC, C. Jordi Girona 1-3, 08034 Barcelona, Spain.\\
$^{3}$ CIMNE, Centre Internacional de M\`etodes Num\`erics a l'Enginyeria, Esteve Terrades 5, E-08860 Castelldefels, Spain.\\
$^*$ Corresponding author.\\
E-mails: {\tt santiago.badia@monash.edu} (SB) 
{\tt pmartorell@cimne.upc.edu} (PM)
{\tt fverdugo@cimne.upc.edu} (FV)
}

\date{\today}

\begin{abstract}
  Unfitted (also known as embedded or immersed) finite element approximations of partial differential equations are very attractive because they have much lower geometrical requirements than standard body-fitted formulations. These schemes do not require body-fitted unstructured mesh generation. In turn, the numerical integration becomes more involved, because one has to compute integrals on portions of cells (only the interior part). In practice, these methods are restricted to level-set (implicit) geometrical representations, which drastically limit their application. Complex geometries in industrial and scientific problems are usually determined by (explicit) boundary representations. In this work, we propose an automatic computational framework for the discretisation of partial differential equations on domains defined by oriented boundary meshes. The geometrical kernel that connects functional and geometry representations generates a two-level integration mesh and a refinement of the boundary mesh that enables the straightforward numerical integration of all the terms in unfitted finite elements. The proposed framework has been applied with success on all analysis-suitable oriented boundary meshes (almost 5,000) in the Thingi10K database and combined with an unfitted finite element formulation to discretise partial differential equations on the corresponding domains.
\end{abstract}

\maketitle

\noindent{{\bf {Keywords}}: Unfitted finite elements, embedded finite elements, clipping algorithms, computational geometry, immersed boundaries, boundary representations.

\section{Introduction}\label{sec:}

Many industrial and scientific applications are modelled by \acp{pde} posed on a non-trivial bounded domain $\Omega$. In these situations, $\Omega$ is described in terms of a boundary representation model (B-REP). These \ac{cad} models are 2-variate, i.e., they are not a parameterisation of $\Omega$ but its boundary $∂ \Omega$; $∂ \Omega$ must be an oriented manifold and $\Omega$ is defined as its interior. On the other hand, the numerical approximation of \acp{pde}, e.g., using \ac{fe} or finite volume schemes, relies on a partition (mesh) of $\Omega$. The traditional simulation pipeline involves unstructured mesh generation algorithms \cite{Si2015,Hu2018}, which take as input the \ac{cad} representation of $∂ \Omega$ and return a mesh covering $\Omega$ (introducing some approximation error). The creation of analysis-suitable \ac{cad} models and body-fitted mesh generation is a non-automatic process that requires intensive human intervention and amounts for most of the simulation time \cite{hughes_isogeometric_2005}. This weak interaction between  (geometry representation) and \ac{cae} functional discretisation is arguably the most serious problem in \ac{cae}, which has motivated the \emph{isogeometric} analysis paradigm \cite{hughes_isogeometric_2005}. Isogeometric analysis has been one of the most active research topic in computational engineering for the last two decades. While this paradigm is sound for \acp{pde} on manifolds (\ac{cad} representations are 2-variate), it does not solve the most ubiquitous situation in practice, i.e., 3D simulations of \acp{pde}  in the $bulk$ of the domain $\Omega$.

Besides, in order to exploit supercomputing resources for unstructured mesh simulations, mesh partitioning strategies must be used, which rely on graph partitioning techniques. Such algorithms are intrinsically sequential and have huge memory requirements \cite{karypis_software_2013}. The mesh partitioning step can easily become the bottleneck (if not a showstopper) of the simulation pipeline for parallel computations on distributed memory machines. Furthermore, the use of such framework in \ac{amr} codes with dynamic load-balancing is not an acceptable option in terms of performance, preventing the use of \ac{amr} in practical large-scale applications with non-trivial geometries. The geometrical discretisation is even more challenging in applications with geometries that evolve in time (like additive manufacturing) or free boundary problems \cite{Neiva2019}, since they require 4D geometrical models (space and time). The generation of body-fitted meshes for complex 3D geometries is still an open problem, and to expect 4D body-fitted generators in a mid term is not reasonable.

In order to solve the current limitations, one could consider unfitted discretisations \cite{mittal_immersed_2005}. A background mesh is used for the discretisation (instead of a body-fitted one) and the geometrical discretisation only requires to generate meshes suitable for integration in the interior defined by $\partial \Omega$  (drastically reducing mesh constraints). An unfitted approach can use tree-based background meshes and exploit scalable and dimension-agnostic mesh generators and partitioners \cite{Badia2020Jun}. Octree-based meshes can be efficiently generated and load-balanced using space-filling curve techniques \cite{bader_space-filling_2012}; see, e.g., the highly scalable $\mathtt{p4est}$ framework \cite{burstedde_p4est_2011} for handling forests of octrees on hundreds of thousands of processors. The extension to a space-time immersed boundaries is feasible since tree-based meshes and marching algorithms are dimension-agnostic.

Unfitted discretisations may lead to unstable and severe ill-conditioned discrete problems~\cite{DePrenter2017} unless a specific technique mitigates the problem. The intersection of a background cell with the physical domain can be arbitrarily small and with unbounded aspect ratio. Despite vast literature on the topic, unfitted finite element formulations that solve these issues are quite recent. Stabilised formulations based on the so-called \emph{ghost penalty} were originally proposed in~\cite{burman2010ghost} for Lagrangian continuous \acp{fe}, and has been widely used since~\cite{burman_cutfem_2015}. The so-called \emph{cell aggregation} or \emph{cell agglomeration} techniques
are an alternative way to ensure robustness with respect to cut location. This approach {is very natural in
\ac{dg} methods, as they can be easily formulated on agglomerated
meshes~\cite{muller2017high}. These techniques have been extended to $\mathcal{C}^{0}$ Lagrangian finite elements in~\cite{Badia2018c} and to mixed methods in \cite{Badia2018a}; the method was coined \ac{agfem}. These unfitted formulations enjoy good numerical properties, such as stability, condition number bounds, optimal convergence, and continuity with respect to data. Distributed implementations for large scale problems have been designed~\cite{Verdugo2019} and  error-driven $h$-adaptivity and parallel tree-based meshes have also been exploited~\cite{Badia2020Jun}. 

Even though unfitted discretisation are motivated by their geometrical flexibility, the current state-of-the-art in unfitted finite elements falls short with respect to the complexity of the geometrics being treated in these publications. The core of the problem is the design of algorithms for the  numerical integration in the interior of background mesh cells only. The vast majority of numerical frameworks rely on implicit level-set descriptions of geometries and marching cubes (or tetrahedra) algorithms, thus limiting their application. We refer the to \cite{Fries2017} for a state-of-the-art review of geometrical discretisation techniques for level-set representations. Geometrical algorithms have also been developed for the intersection of 3D tetrahedral meshes for the unfitted discretisation of interface problems (see \cite{Massing2013,Johansson2019} and references therein).  

The main motivation of this work is to provide a new geometrical framework that covers all the needs of  unfitted techniques and is amenable to arbitrarily complex 3D geometries {represented by \ac{stl} meshes, i.e., oriented faceted linear surface representations. This is one of the most common situations in \ac{cae}, in which \ac{stl} meshes are used to define complex objects.}  In particular, the starting point of the algorithm, as in unstructured mesh generation, is a boundary mesh for $\partial \Omega$. We design an algorithm for computing the intersection of each cell in a background mesh and the interior of the boundary mesh. The number of faces in the boundary mesh intersecting a background cell can be in the order of hundreds or even thousands for very complex geometrical representations. As a result, the proposed algorithm must be resilient to rounding errors and provide answers accurate up to machine precision in all cases. With this algorithm, we complete an automatic simulation framework that takes a standard \ac{cad} representation, {an \ac{stl} mesh,} and returns the \ac{pde} solution obtained from an unfitted discretisation. The procedure is fully automatic and allows us to exploit all the benefits of unfitted formulations described above on complex geometries defined by {\ac{stl}} representations. On the other hand, with the proposed formulation, the geometrical error (determined by the boundary mesh) and the functional error (determined by the background mesh) are completely decoupled. This remarkable property is not shared by \acp{fe} on unstructured meshes (both geometry and functional discretisation rely on the same geometrical discretisation) or standard level-set approaches with marching algorithms on background cells (the geometrical approximation is determined by the background mesh). {We note that the geometrical framework proposed in this work can readily be applied to other numerical techniques that can be posed on general polytopal meshes, e.g., hybridised formulations on agglomerated meshes \cite{2109.09983} or mollified \acp{fe} \cite{Febrianto2021}.}

The first key ingredient of the proposed framework is a robust clipping algorithm for convex polytopes. A popular method for clipping is the one by Sutherland and Hodgman in \cite{Sutherland1974} (see also~\cite{Stephenson1975Sep}). Recent implementations and improvements of these algorithms can be found in~\cite{Lpez2018} and an extension to non-convex polyhedra can be found in~\cite{Lpez2019}. These methods require accuracy checks and the handling of all degenerate branches. Instead, the approach proposed by Sugihara and co-workers, called \emph{combinatorial abstraction}, is an example of a numerically robust scheme for the intersection of convex polyhedra~\cite{Sugihara1994}. An implementation of this algorithm has been proposed in~\cite{Powell2015} for intersecting a tetrahedron and a background Cartesian mesh. Still, Sugihara's method relies on assumptions that are not true in general and current implementations are not designed to deal with a large number of clipping planes or non-convex geometries. In this work, we build on \cite{Sugihara1994,Powell2015} to design a robust  algorithm and an efficient implementation for the clipping of a polyhedron and a plane that can naturally handle with degenerate possibly non-connected and non-convex outputs and is suitable for our specific target. 

Since the boundary representation is not convex in general, the second key ingredient is a convex decomposition algorithm that transform a non-convex intersection into a set of convex ones. Different algorithms have been proposed for the convex decomposition of non-convex polyhedra~\cite{Chazelle1984,Hachenberger2008}. In this work, we use polyhedron decomposition ideas for the intersection problem at hand. The original intersection problem is decomposed into a set of convex clipping problems for which one can use our convex clipping strategy.

Finally, it is essential to design mechanisms that provide robustness of the algorithm with respect to rounding errors. The main problem is the potentially huge number of clipping planes to be processed. Our methods are based on a discrete level-set representation of planes (instead of a more traditional parametric representation) and a specifically oriented graph representations of polyhedra. The motivation for this choice is to maximise symbolic computations and define geometrical operations that are numerically robust under rounding errors. We also provide techniques that identify and merge quasi-aligned planes. 

The outcomes of this article are the following:
\begin{itemize}
	\item A robust and efficient intersection algorithm for computing the interior of cells given a boundary mesh representation;
	\item The combination of the intersection algorithm with unfitted \acp{fem} for a body-fitted mesh free computational framework that is applicable to the discretisation of \acp{pde} on explicit representations of complex geometries;
	\item A detailed robustness analysis of the geometrical algorithms on the Thingi10K database with about 5,000 surface meshes \cite{Zhou2016};
	\item The numerical experimentation of an unfitted \ac{fe} solver that relies on the proposed geometrical intersection engine;
	\item A performance analysis of the proposed framework and an open-source implementation~\cite{Martorell_STLCutters_2021}.
\end{itemize}

The outline of the article is as follows.  In Section~\ref{sec:fem}, we present some unfitted \ac{fe} methods and their geometrical requirements. In Section~\ref{sec:intersection}, we provide the geometrical algorithm that computes the intersection of background cells and oriented surface meshes. In Section~\ref{sec:num-exp}, we report a thorough numerical experimentation of the proposed algorithms on almost 5,000 meshes in the Thingi10K database~\cite{Zhou2016}. We show the remarkable robustness of the geometrical algorithm, providing very low geometrical error quantities in all cases. The algorithms are combined with unfitted \ac{fe} methods to approximate \acp{pde} on these complex geometries, showing the expected convergence orders of accuracy. We also analyse the computational performance of the proposed framework and provide details about the corresponding open source implementation~\cite{Martorell_STLCutters_2021}. Finally, some conclusions and future work lines are drawn in Section~\ref{sec:conclusions}.

\section{Unfitted finite element discretisations}\label{sec:fem}

Let us consider an open Lipschitz domain $\Omega \subset \mathbb{R}^3$ (the 2D case is an obvious restriction) in which we want to approximate a system of \acp{pde}.   In this work, we are interested in domains that are described as the interior of an oriented surface polygonal mesh $\mathcal{B}$ of $\partial \Omega$. \acp{pde} usually involve Dirichlet boundary conditions on $\Gamma_D$ and Neumann boundary conditions on $\Gamma_N$, where $\Gamma_D$ and $\Gamma_N$ are a partition of $\partial \Omega$. Such partition must be respected by the geometrical representation, e.g., the {\ac{stl}} model. Thus, we consider that $\mathcal{B}_D$ and $\mathcal{B}_N$ are geometrical discretisations of $\Gamma_D$ and $\Gamma_N$, resp., and $\mathcal{B} \doteq \mathcal{B}_D \cup \mathcal{B}_N$.  

Our motivation in this work is to enable the use of grid-based unfitted numerical schemes that are automatically generated from $\mathcal{B}$, due to its industrial and scientific relevance. Embedded discretisation techniques alleviate geometrical constraints, because they do not rely on body-fitted meshes. Instead, these techniques make use of a background partition $\mathcal{T}^{\mathrm{bg}}$ of an arbitrary artificial domain $\Omega^{\mathrm{art}}$ such that $\Omega \subset \Omega^{\mathrm{art}}$. The artificial domain can be trivial, e.g., it can be a bounding box of $\Omega$. Thus, the computation of $\mathcal{T}^{\mathrm{bg}}$ is much simpler (and cheaper) than a body-fitted partition of $\Omega$. In this work, we consider a Cartesian mesh $\mathcal{T}^{\mathrm{bg}}$ for simplicity in the exposition, even though the proposed approach could readily be extended, e.g., to a tetrahedral structured background mesh obtained after simplex decomposition. 

The abstract exposition of unfitted formulations considered in this work is general and accommodates different unfitted \ac{fe}  techniques that have been proposed in the literature, e.g.~the \ac{xfem}~\cite{belytschko_arbitrary_2001} (for unfitted interface problems), the cutFEM method~\cite{burman_cutfem_2015} based on ghost penalty stabilisation, the \ac{agfem}~\cite{Badia2018c}, the finite cell method~\cite{Schillinger2015} and \ac{dg} methods with cell aggregation~\cite{muller2017high}, to mention a few. 

In order to define a \ac{fe} space on unfitted meshes, we do the following cell classification. The cells in the background partition with null intersection with $\Omega$ are \emph{exterior} cells. The set of exterior cells in $\mathcal{T}^{\mathrm{bg}}$ is denoted by $\mathcal{T}^{\mathrm{out}}$  is not considered in the functional discretisation and can be discarded.  $\mathcal{T} \doteq \mathcal{T}^{\mathrm{bg}} \setminus \mathcal{T}^{\mathrm{out}}$ is the \emph{active} mesh (see \fig{fig:immersed-setup-a}). The above mentioned techniques make use of standard \ac{fe} spaces on $\mathcal{T}$ to define the finite-dimensional space $V$ in which to seek the solution and also test the weak form of the \ac{pde}. The unfitted problem reads as follows: find $u \in V$ such that 
\begin{equation}
	a(u,v) = \ell(v), \qquad \forall v \in \mathcal{V}, 
\end{equation} 
where 
\begin{equation}
a(u,v) = \int_{\Omega}  L_\Omega(u,v) \mathrm{d} \Omega + \int_{\Gamma^{\mathrm{D}}}^{} L_{D} (u,v) \mathrm{d} \Gamma + \int_{\mathcal{F}}^{} L_{\mathrm{sk}}(u,v) \mathrm{d} \Gamma,
\end{equation}
and
\begin{equation}
\ell(v) =  \int_{\Omega}  F_\Omega(v) \mathrm{d} \Omega + \int_{\Gamma^{\mathrm{N}}} F_N (v) \mathrm{d} \Gamma + \int_{\Gamma_D} F_{D} (v) \mathrm{d} \Gamma.
\end{equation}
The bulk terms $L_\Omega$ and $F_\Omega$ include the differential operator (in weak sense), the source term and possibly some other numerical stabilisation terms. The operators $L_D$ and $F_D$ integrated on $\Gamma_D$ represent the terms related to the weak imposition of Dirichlet boundary conditions, e.g., the so-called Nitsche's method, which is commonly used in unfitted formulations. The term $F_N$ on $\Gamma_N$ represents the Neumann boundary conditions of the problem at hand. We denote with $\mathcal{F}$ the skeleton of the active mesh, i.e., the set of interior faces of $\mathcal{T}$. The term $L_{\mathrm{sk}}$  collects additional penalty terms that include weak imposition of continuity in \ac{dg} methods or ghost penalty stabilisation techniques.

\definecolor{pacyan}{RGB}{104,207,207}
\definecolor{pablue}{RGB}{0,68,170}
\definecolor{payellow}{RGB}{231,211,68}
\definecolor{pared}{RGB}{255,0,0}

\begin{figure}[http]
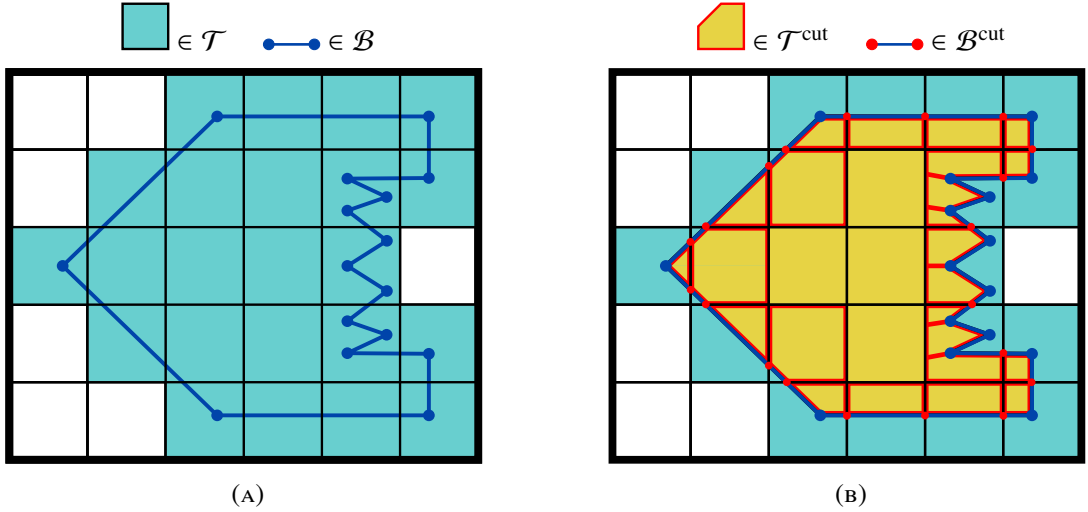

  \begin{subfigure}[b]{0.49\textwidth}
    \centering
    \tikz{ \draw[draw=black,fill=pacyan, line width=0.9pt] (0,0) rectangle (0.6,0.6);   } $\in\mathcal{T}$ \quad \tikz{ \draw[line width=1pt, color=pablue] (0,0) -- (0.6,0.0); \draw[fill,color=pablue] (0,0) circle (0.07);  \draw[fill,color=pablue] (0.6,0.0) circle (0.07);  } $\in\mathcal{B}$\\[0.5em]    
    \includefig{0.8\textwidth}{T_2d}
    \caption{}
    \label{fig:immersed-setup-a}
  \end{subfigure}
  \begin{subfigure}[b]{0.49\textwidth} 
    \centering
    \tikz{ \draw[draw=pared,fill=payellow, line width=0.9pt] (0,0) -- (0.6,0.0) -- (0.6,0.6) -- (0.3,0.6) -- (0.0,0.3) -- cycle ;   } $\in\mathcal{T}^{\mathrm{cut}}$ \quad \tikz{ \draw[line width=1pt,color=pablue] (0,0) -- (0.6,0.0); \draw[fill,color=pared] (0,0) circle (0.07);  \draw[fill,color=pared] (0.6,0.0) circle (0.07);  } $\in\mathcal{B}^{\mathrm{cut}}$\\[0.5em]   
    \includefig{0.8\textwidth}{Tcut_2d}
    \caption{}
  \end{subfigure}
  \caption{Example of an embedded non-convex domain in 2D. Left hand side figure (a) shows an active mesh $\mathcal{T}$ and an oriented skin mesh $\mathcal{B}$. The two-level meshes in (b), namely $\mathcal{T}^{\mathrm{cut}}$ and $\mathcal{B}^{\mathrm{cut}}$, are computed using the techniques proposed in this work to integrate unfitted formulations. Cut cells in $\mathcal{T}$ are split into a set of convex polytopes in $\mathcal{T}^{\mathrm{cut}}$. We note that $\mathcal{T}^{\mathrm{cut}}$ is not conforming across cells in 3D in general (only in 2D).}
  \label{fig:immersed-setup-b}
\end{figure}
    
Since \ac{fe} methods are piecewise polynomials, the integration of all this terms rely on a cell-wise decomposition (of bulk and surface terms). However, in order to respect the geometry and solve the \ac{pde} on the right domain, one must perform these integrals on domain interiors. In particular, we have 
\begin{equation}\label{eq:dom-terms}
\int_{\Omega} (\cdot) \mathrm{d} \Omega = \sum_{K \in \mathcal{T}}  \int_{K \cap \Omega }^{} (\cdot) \mathrm{d} \Omega.
\end{equation}

As commented above, the surface mesh $\mathcal{B}$ in which we aim to integrate the boundary terms and the background mesh that defines the cell-wise polynomial \ac{fe} functions are not connected, i.e., one is not the boundary restriction of the other. Thus, the integration of boundary terms must be computed cell-wise as follows: 
\begin{equation}\label{eq:bou-terms}
	\int_{\Gamma_*} (\cdot) \mathrm{d} \Gamma = \sum_{K \in \mathcal{T}} \sum_{F \in \mathcal{B}_*} \int_{F \cap K}^{} (\cdot) \mathrm{d} \Gamma, \qquad * \in \{ D, N \}.
\end{equation}
We note that, even though it does not represent any problem for the machinery we want to propose here, skeleton terms (common in ghost penalty and \ac{dg} methods) can still be integrated on the whole skeleton faces, and there is no need in general to reduce these integrals to the domain interior. 

As a result, the \emph{only} geometrical complication of unfitted finite element schemes is the integral over $K \cap \Omega$ and $K \cap \mathcal{B}$ (or more specifically, $\mathcal{B}_N$ and $\mathcal{B}_D$). Such operations (and specially the first one) are hard for $\partial \Omega$ representing a complex shape explicitly determined by an {\ac{stl}} model and will be the target of the next section.
Let us stress the fact that the tools described below are \emph{only used for integration purposes}. The functional spaces are defined in the background mesh, which provides lots of flexibility (no inter-cell consistency or shape regularity requirements) compared to unstructured mesh generation.

In particular, the geometrical queries that are required by unfitted formulations can be solved as follows. First, the intersection of the surface cells in $\mathcal{B}$ against all background cells in $\mathcal{T}$ produces a new surface mesh $\mathcal{B}^\mathrm{cut}$ that is a refinement of $\mathcal{B}$ that describes the same geometry (up to machine precision). It can be indexed as a two-level mesh, in which first one computes its portion for a cell $K \in \mathcal{T}$, $\mathcal{B}^{\mathrm{cut}}_K \doteq \{ S \cap K  : S \in \mathcal{B} \}$
and $\mathcal{B}^\mathrm{cut} \doteq \bigcup_{K \in \mathcal{T}} \mathcal{B}^{\mathrm{cut}}_K$. Thus, $\mathcal{B}^\mathrm{cut}$ can readily be used for the integration of the boundary terms in (\ref{eq:bou-terms}); we can analogously use $\mathcal{B}_D$ (resp., $\mathcal{B}_N$) to produce $\mathcal{B}^{\mathrm{cut}}_D$ (resp., $\mathcal{B}^{\mathrm{cut}}_N$). In any case, $\mathcal{B}^\mathrm{cut}$ preserves global conformity in 3D and can also be understood as a standard one-level polytopal mesh. Second, for each cell $K \in \mathcal{T}$, we want to compute a sub-mesh $\mathcal{T}_K$ (composed of convex polyhedra) of the interior of the cell, i.e., $K \cap \Omega$. We represent with $\mathcal{T}^{\mathrm{cut}} \doteq \{ \mathcal{T}_K : K \in \mathcal{T} \}$ the resulting two-level mesh such that $\bigcup_{K \in \mathcal{T}} \bigcup_{L \in \mathcal{T}_K} L = \Omega$ (up to machine precision). This two-level \emph{integration} mesh can readily be used to compute the integrals in (\ref{eq:dom-terms}). $\mathcal{T}^{\mathrm{cut}}$ is not conforming across background cells, since this mesh is is only needed for the numerical integration. Since the result of this algorithm is a set of meshes $\mathcal{T}_K$ composed of general convex polytopes, we can now use numerical quadratures for the integration on these polytopes. For these purposes, one can use quadrature rules on general polytopes (see, e.g., \cite{Chin2020}) or a straightforward simplex decomposition and standard quadrature rules on triangles/tetrahedra. \fig{fig:immersed-setup-b} illustrates the  construction of $\mathcal{B}^{\mathrm{cut}}$ and $\mathcal{T}^{\mathrm{cut}}$.
 
\section{Intersection algorithms}\label{sec:intersection}

In this section, we provide an algorithm that given the background mesh $\mathcal{T}$ and the oriented surface polygonal mesh $\mathcal{B}$ (resp., $\mathcal{B}_D$ and $\mathcal{B}_{N}$), it returns $\mathcal{T}^{\mathrm{cut}}$ and $\mathcal{B}^\mathrm{cut}$ (resp., $\mathcal{B}^\mathrm{cut}_D$ and $\mathcal{B}^\mathrm{cut}_N$). The problem when computing these meshes is the fact that the geometrical intersection and parametric distance computation algorithms are in general not robust for these purposes, due to inexact arithmetic. In this work, we aim at designing an algorithm that is robust and can be readily applied to any surface mesh with a well-defined interior. In order to attain such level of robustness, we work at different levels:

\begin{itemize}

  \item First, we provide in Section~\ref{sec:definitions} and Section~\ref{sec:half-space} a computer representation of polyhedra (or surfaces) and planes, resp., that are suitable for intersection algorithms with inexact arithmetic.
  
  \item Second, in Section~\ref{sec:poly-plane-intersection} we discuss a novel algorithm for the intersection of polyhedra and half-spaces. The algorithm is inspired by Sugihara's intersection algorithm~\cite{Sugihara1994} and Powell and Abell implementation in \cite{Powell2015}, but departs from these algorithms to make it suitable for our specific purposes. We provide the complete algorithm (up to minor implementation details) for the intersection of a convex polyhedron and a plane. 
  
  \item Third, in Section~\ref{sec:intersection-P-S} we consider the intersection of a convex polyhedron against a \emph{non-convex} surface. In order to do that, we need to define a recursive convex decomposition algorithm that re-states the original intersection problem as a set of intersections between convex polyhedra and surfaces, for which we can use the algorithms in Section~\ref{sec:poly-plane-intersection}. 
  
  \item Fourth, in order to have a robust algorithm, it is not enough with the proposed representation of polyhedra and planes and the proposed intersection algorithms. A common problem that appears in inexact arithmetic is the case of multiple planes that are quasi-aligned (conceptually, aligned up to machine precision). It has been proven that merging (enforcing the planes to be exactly aligned) dramatically improves the robustness of the overall algorithm. The approach we propose to merge planes is presented in Section~\ref{sec:quasi-coplanarity}.
  
  \item Finally, with all these ingredients, we can design the global intersection algorithm in Section~\ref{sec:global}, which returns $\mathcal{B}^{\mathrm{cut}}$ and $\mathcal{T}^{\mathrm{cut}}$ explained above. 

\end{itemize}

\subsection{Polyhedra and polygonal surface representations}\label{sec:definitions}

In this work, we have to deal with hundreds (or even thousands in some limit cases) of planes clipping a cell. This situation makes robustness essential, which prevents us from using methods that require accuracy checks and the handling of all degenerate branches. For this reason, our starting point is Sugihara's method \cite{Sugihara1994} and its implementation in~\cite{Powell2015}. However, the objective in~\cite{Powell2015} is to intersect a tetrahedron and a Cartesian mesh and thus restricted to a convex surface mesh and polyhedron. We propose below an algorithm that keeps robustness for a large number of clipping planes and non-convex situations.

Sugihara's method relies on the following assumption. A convex polyhedron intersected by a plane must produce two connected polyhedrons. As Sugihara pointed out in his seminal work, this is not the case in inexact arithmetic. In order to expose the problem, let us consider a polyhedron face with more than three vertices. In exact arithmetic, all these points belong to the same plane. However, this is not true in numerical computations, \emph{co-planarity is only true up to machine precision}. (Below, we use the prefix \emph{quasi-} to indicate a geometrical concept that is true for exact arithmetic but only approximate in finite precision.) Next, let us consider an oriented plane that is quasi-coplanar to the face up to machine precision. The classification of a vertex as interior, exterior, or on the plane is completely determined by rounding errors, thus unreliable. E.g., it can lead to a non-connected partition of the polyhedron that is impossible in exact arithmetic. In singular cases, Sugihara proposes an algorithm to re-classify the vertices on the two sides of the cutting plane based on logical arguments.

Despite Sugihara's method, we do not consider any re-classification of vertices to satisfy Sugihara's assumption in \cite{Sugihara1994}. We consider an algorithm for the clipping of a polyhedron and a plane that can naturally handle possibly non-connected and non-convex outputs. In any case, the loss of convexity can only produce rounding errors and the resulting polyhedron is \emph{quasi-convex}.

Let us start introducing some basic notation about graphs. Given a graph $G$, we denote with $\mathtt{vert}(G)$ the set of vertices of the graph and with $\mathtt{adj}(G)$ the adjacencies. The adjacency of a vertex $\alpha \in \mathtt{vert}(G)$, i.e., the set of vertices connected to $\alpha$ by an edge of the graph $G$, is denoted by $\mathtt{adj}(G)(\alpha)$. We can extract the set of \emph{ connected components} (or sub-graphs) $\mathtt{comp}(G)$ of a graph $G$. 

A graph can readily be constructed from a set of vertices $V$ and  the vertices adjacencies $E$; we represent this construction with $\mathtt{graph}(V,E)$. We can also make use of a constructor $\mathtt{graph}(V,\mathtt{C})$, where $\mathtt{C}$ is a condition that determines whether two vertices $\alpha, \beta \in V$ are connected ($\mathtt{C}(\alpha,\beta)$ is true) or not. This construction allows us to define both directed and undirected graphs, for symmetric or non-symmetric conditions, respectively.

In order to represent polyhedra, we need to make use of \emph{rotation systems}, i.e., a sub-type of graphs in which the adjacency of each vertex is a \emph{cyclic order}. In a rotation system $R$, given a vertex $\alpha \in \mathtt{vert}(R)$ and $\beta \in \mathtt{adj}(R)(\alpha)$, there is a well-defined \emph{previous} and $next$ in $\mathtt{adj}(R)(\alpha)$, defined by the cyclic ordering. Thus, we can define $\mathtt{next}(\alpha;\beta)$ as the vertex after $\beta$ in the cyclic ordering $\mathtt{adj}(R)(\alpha)$.

\begin{definition}[Polyhedron representation]\label{def:pol}
The boundary of a polyhedron $P$ is an oriented closed surface made of polygons in which the edges around a vertex admit a cyclic ordering that encodes the surface orientation. The cyclic ordering of the adjacency (neighbours) of a vertex $\alpha \in \mathtt{vert}(P)$ is determined by the clockwise ordering of edges as observed when positioned outside of $P$ on $\alpha$. Thus, a polyhedron $P$ can be represented as a rotation system whose vertices are points in $\mathbb{R}^3$. This description of a polyhedron has been exploited in \cite{Powell2015}. Figure~\ref{fig:alg2-simple-intersection} at step $(i)$ shows a cube representation as a rotation system with clockwise ordering of neighbours.
\end{definition}

We can define a specific traversal of the polyhedron vertices using the definition of $\mathtt{next}$ defined by the cyclic ordering above. Given an edge $(\alpha_0,\alpha_1)$ of the polyhedron, subsequent vertices repeatedly applying $\alpha_{i+1} \gets \mathtt{next}(\alpha_i;\alpha_{i-1})$. The faces of the polyhedron are the \emph{closed paths} determined by this graph traversal, i.e., a face is defined by $\alpha_0, \alpha_1$ and the iteration $\alpha_{i+1} \gets \mathtt{next}(\alpha_i;\alpha_{i-1})$ till the result is $\alpha_0$; the face is a 2D polygon itself. We represent the set of faces in a polyhedron with $\mathtt{faces}(P)$. 

An open polygonal oriented surface $\vec{\Gamma}$ can also be represented as the polyhedron plus information about which vertices lie on the boundary, which are represented with $\mathtt{bou}(\vec{\Gamma})$. It is convenient to \emph{close} these open surfaces. We define the concept of \emph{open} vertex $o$. Conceptually, $o$ is a vertex at infinite distance of the surface and exterior to all the faces of the surface mesh $\vec{\Gamma}$. Algorithm~\ref{alg:surface-to-polyhedron} receives a surface mesh $\vec{\Gamma}$ and returns a polyhedron by modifying the surface graph by appending the artificial \emph{open} node to the adjacency of boundary vertices (line~\ref{ln:pol-1}). As we want closed paths to represent polyhedron faces, we need a mechanism to avoid  vertices in $\mathtt{bou}(\vec{\Gamma})$ to define a closed path; \emph{open} vertices break this path. This construction is illustrated in Figure~\ref{fig:alg1-open-close}.

\begin{algorithm}[H]
	\caption{$\mathtt{pol}(\vec{\Gamma})$}
	\justifying
	\begin{algorithmic}[1]
    \STATE ${V} \gets \mathtt{vert}(\vec{\Gamma}), \quad {E} \gets \mathtt{adj}(\vec{\Gamma}), \quad  \partial {V} \leftarrow \mathtt{bou}(\vec{\Gamma}), \quad {V} \leftarrow {V} \cup \{o\}$
	\FOR{$v \in \partial {V}$}
	    \STATE ${E}(v) \leftarrow ({E}(v),o) $ \label{ln:pol-1}
	\ENDFOR
	\RETURN $\mathtt{graph}({V},{E})$ 
	\end{algorithmic}\label{alg:surface-to-polyhedron}
  \end{algorithm}

\begin{figure}[http]
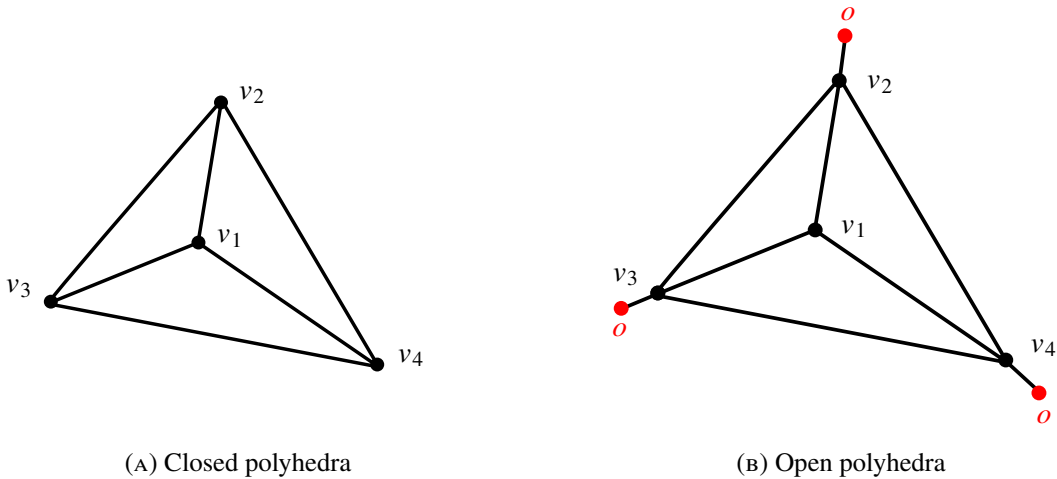

  \begin{subfigure}[b]{0.49\textwidth}
    \centering
    \includefig[\normalsize]{0.8\textwidth}{alg1_a}
    \caption{Closed polyhedra}
  \end{subfigure}
  \begin{subfigure}[b]{0.49\textwidth}
    \centering
    \includefig[\normalsize]{0.8\textwidth}{alg1_b}
    \caption{Open polyhedra}
  \end{subfigure}
  \caption{Example of \alg{alg:surface-to-polyhedron} that converts (a) a closed polyhedra into (b) an open surface polyhedra by adding \emph{open} nodes to the boundary vertices $\partial V = \left\{v2,v3,v4\right\}$. The boundary vertices represent a graph cycle in (a), i.e., a polyhedron face. However, that cycle is broken in (b) because of the edges to \emph{open} node, $o$. Hence, the polyhedron is open in (b) and represents a surface, as intended.}
  \label{fig:alg1-open-close}
\end{figure}
\subsection{Half-space representations}\label{sec:half-space}

Given an oriented plane $\vec{\pi}$, i.e., determined by a face of the polyhedron, one can define its corresponding \emph{open} half-space as the set of points in the interior side of the plane. We use the following discrete representation of this space for a given set of vertices.

\begin{definition}[Half-space representation]
  We represent the half-space corresponding to an oriented plane $\vec{\pi}$ using a \emph{discrete level-set} $\boldsymbol{h}$ with respect to a set of vertices $V$, i.e., the set of signed distances of the vertices in $V$ to the plane $\vec{\pi}$. $\boldsymbol{h}$ can be represented as an array of real numbers of length $|V|$ . We use the convention that a negative sign means interior point (positive for exterior points). We assume that half-spaces are open and define their closure as $\overline{\boldsymbol{h}}$.  We note that the only difference between the open and closed half-spaces is the definition of $\in$; vertices with zero distance belong to the closed half-space but not to the open one. Let us represent the plane with inverted orientation as $-\vec{\pi}$. Given the half-space $\boldsymbol{h}$ of $\vec{\pi}$, we define the one for $-\vec{\pi}$ as $-\boldsymbol{h}$. The complement of $\boldsymbol{h}$ is $\overline{-\boldsymbol{h}}$, i.e., $\boldsymbol{h} \oplus \overline{-\boldsymbol{h}} = \mathbb{R}^3$. 
\end{definition}

{The discrete level-set representation of a plane with respect to a set of vertices $V$ can be determined by computing all the signed distances between vertices in $V$ and the plane. Let us consider a set $S$ of planes. We represent the corresponding set of half-spaces} as a \emph{signed distance matrix} $\mathbf{H}$ in which the rows are half-spaces (the first index is the plane in $S$) and the columns are signed distances to all planes for a given vertex in $V$ (the second index is the vertex). We also need to use block partitions of the matrix. E.g, if $S = X \cup Y$ and $V = W \cup Q$, we use the notation $\mathbf{H}_{XY,WQ}$ for the whole matrix while the blocks are represented using specific subscripts, e.g., the matrix block for planes $X$ and vertices $W$ is $\mathbf{H}_{X,W}$. We abuse of notation when dealing with polyhedra and distance matrices. E.g., given a polyhedron $S$, we use $\mathbf{H}_{{S},*}$ instead of $\mathbf{H}_{\mathtt{faces}(S),*}$ (or more accurately, the planes that contain $\mathtt{faces}(S)$) and  $\mathbf{H}_{*,S}$ instead of $\mathbf{H}_{*,\mathtt{vert}(S)}$; the symbol $*$ means all indices in that dimension.

\subsection{Clipping a polyhedron with a plane}\label{sec:poly-plane-intersection}

We are in position to provide Algorithm~\ref{alg:polyhedron-intersection}, in which we compute the clipping of a convex polytope $P$ against a half-space $\boldsymbol{h} \in \rm{rows}(\mathbf{H})$. 
The result of this algorithm is (i) the new polytope obtained after clipping $P$ with $\boldsymbol{h}$ and (ii) the new set of half-spaces $\mathbf{H}$  after eliminating $\boldsymbol{h}$, eliminating distances to vertices that are not in $P$ anymore and adding distances to newly created vertices. The half-space $\boldsymbol{h}$ can be open or closed; it does not affect $\mathbf{H}$, since the same distances are required in both cases. {The main steps in Algorithm~\ref{alg:polyhedron-intersection} are illustrated in Figure~\ref{fig:alg2-simple-intersection}, in which a cube is clipped by a plane.}

Let us assume that $\boldsymbol{h}$ is open. The algorithm iterates over interior vertices (lines~\ref{ln:polyhedron-intersectionint-2}-\ref{ln:polyhedron-intersectionint-3}). At each vertex $\alpha \in \boldsymbol{h}$, we look for vertices in the adjacency of $\alpha$ that are exterior (\ref{ln:polyhedron-intersectionint-4}-\ref{ln:polyhedron-intersectionint-5}).\footnote{In line~\ref{ln:polyhedron-intersectionint-4}, we use the notation $\mathtt{enum}$ over an iterator to describe a new iterator that yields a tuple $(i,a)$ in which $i$ is a counter starting at 1 and $a$ is the $i$-th value from the given iterator.} If we find an exterior vertex $\beta$ in the adjacency of $\alpha$. We have found an edge $(\alpha,\beta)$ that intersects the plane related to the half-space. In line~\ref{ln:polyhedron-intersectionint-6}, we compute the coordinates of the intersection vertex $\delta$ computing its coordinates as:
\begin{equation}\label{eq:edge-plane-vertex}
\boldsymbol{x}_\delta 
\doteq \xi_{\alpha \delta} \boldsymbol{x}_{\alpha} + \xi_{\beta\delta} \boldsymbol{x}_{\beta}
\doteq \frac{-\boldsymbol{h}_\beta}{\boldsymbol{h}_\alpha-\boldsymbol{h}_\beta }  \boldsymbol{x}_\alpha
+
 \frac{\boldsymbol{h}_\alpha}{\boldsymbol{h}_\alpha-\boldsymbol{h}_\beta } \boldsymbol{x}_\beta.
\end{equation}
This computation has also been used in \cite{Powell2015} because it is much more robust than using an intersection algorithm between planes and edges in parametric form. The denominator is always positive, $\xi_{\alpha\delta}$ and $\xi_{\beta\delta}$ are non-negative by construction and the resulting vertex $\delta$ always lies between $\alpha$ and $\beta$. Furthermore, the computation of the signed distance between the new vertex $\delta$ and any half-space $\boldsymbol{h}' \in \mathbf{H}$ can readily be computed as:
\begin{equation}\label{eq:dist-delta-plane}
  \boldsymbol{h}'_\delta = \xi_{\alpha \delta} \boldsymbol{h}'_\alpha + \xi_{\beta\delta} \boldsymbol{h}'_\beta.
\end{equation} 
In order to illustrate the robustness of this approach, let us discuss what happens in the singular case in which a vertex is exactly on the intersecting plane. Using the definition above for an open half-space, this vertex is exterior. Any edge that connects it to an interior node will be intersected and a new vertex will be inserted. The new vertex distance to the planes (and coordinates) will have \emph{exactly} the same coordinates as the vertex on the boundary, due to the expression in (\ref{eq:edge-plane-vertex}). In line~\ref{ln:polyhedron-intersectionint-7} we create the adjacency (cyclic order) for $\delta$ with $\alpha$ in the first position and two additional positions not defined yet. On the other side, we replace the intersected edge $(\alpha,\beta)$ with $(\alpha,\delta)$.

After this loop, we have identified all intersected edges, computed the new vertices after the intersection and modified the adjacencies. The adjacencies of the new vertices are not yet complete because we have not included the edges on the new face created after clipping; the edges of this face only include new vertices. We perform a new loop over new vertices in line~\ref{ln:polyhedron-intersectionint-12}. For each new vertex $\alpha$,  we start a graph traversal (lines~\ref{ln:polyhedron-intersectionint-13}-\ref{ln:polyhedron-intersectionint-16}) till we find another new vertex $\beta$. Then, we put $\beta$ in the second position of $\mathtt{adj}(P)(\alpha)$ and $\alpha$ in the third position of $\mathtt{adj}(P)(\beta)$. When this loop finishes, we have the complete adjacencies of the new vertices.     

It only remains to add to the polytope the new vertices (and their adjacencies) and to eliminate the exterior vertices (line~\ref{ln:polyhedron-intersectionint-19}-\ref{ln:polyhedron-intersectionint-20}). We also compute the signed distances of new vertices to half-spaces in $\mathbf{H}$ and eliminate the ones related to exterior vertices (inserting/removing rows to this matrix in line~\ref{ln:polyhedron-intersectionint-21}). The distance to any half-space are computed using (\ref{eq:dist-delta-plane}). Since $P$ has been clipped with $\boldsymbol{h}$, it is eliminated from $\mathbf{H}$.   

The most salient property of this algorithm is that most computations are symbolic, with the only exception of the new vertex coordinates in line~\ref{ln:polyhedron-intersectionint-6} and new distances in line~\ref{ln:polyhedron-intersectionint-21}. However, the computation of these quantities has already been designed in such a way that they are well-posed in finite precision, using the expressions (\ref{eq:edge-plane-vertex})-(\ref{eq:dist-delta-plane}) discussed above.

We have considered the intersection with one half-space. But we can recursively use the algorithm to intersect with multiple planes, since we do not assume any specific topology of the initial quasi-convex polyhedron. With minor modifications, one can also extract not only the interior but also the exterior graph at the same time, reusing computations.

\begin{algorithm}[h]
	\caption{$({P},\mathbf{H}) \cap \boldsymbol{h}$}
	\justifying
	\begin{algorithmic}[1]
	\STATE $V_P \gets \mathtt{vert}(P), \quad E_P \gets \mathtt{adj}(P); \quad {V}_P^{\mathrm{new}} \leftarrow \emptyset$ \label{ln:polyhedron-intersectionint-1}
	\FOR{$\alpha \in V_P$} \label{ln:polyhedron-intersectionint-2}
	\IF{$\alpha \in \boldsymbol{h}$} \label{ln:polyhedron-intersectionint-3}
	\FOR{$(i,\beta) \in \mathtt{enum}(E_P(\alpha))$} \label{ln:polyhedron-intersectionint-4}
	  \IF{$\beta \not\in \boldsymbol{h}$}  \label{ln:polyhedron-intersectionint-5}
	    \STATE $\delta \leftarrow \vec{\alpha \beta} \cap \boldsymbol{h}$ 
	    \label{ln:polyhedron-intersectionint-6}
	    \STATE ${V}_P^{\mathrm{new}} \leftarrow V_P^{\mathrm{new}}\cup \{ \delta \}$; \quad $E_P \leftarrow E_P\cup \{ (\delta , (\alpha,\varnothing,\varnothing)) \}$; \quad  $E_P(\alpha)[i] \leftarrow  \delta$ \label{ln:polyhedron-intersectionint-7}
	  \ENDIF \label{ln:polyhedron-intersectionint-8}
	  \ENDFOR \label{ln:polyhedron-intersectionint-9}
	  \ENDIF \label{ln:polyhedron-intersectionint-10}
	\ENDFOR \label{ln:polyhedron-intersectionint-11}
	\FOR{$\alpha \in V_P^{\mathrm{new}}$} \label{ln:polyhedron-intersectionint-12}
	  \STATE $(\beta,\delta) \leftarrow (\alpha,E_P(\alpha)[1])$ \label{ln:polyhedron-intersectionint-13}
	  \WHILE{$\delta  \not\in V_P^{\mathrm{new}}$} \label{ln:polyhedron-intersectionint-14}
	  \STATE $(\beta,\delta) \leftarrow (\delta,\mathtt{next}(\delta;\beta))$ \label{ln:polyhedron-intersectionint-15}
	  \ENDWHILE  \label{ln:polyhedron-intersectionint-16}
	  \STATE $E_P(\alpha)[2] \leftarrow \delta$, \quad $E_P(\delta)[3] \leftarrow  \alpha$ \label{ln:polyhedron-intersectionint-17}
	\ENDFOR \label{ln:polyhedron-intersectionint-18}
	\STATE ${V}^{\mathrm{out}} \leftarrow \left\{ \alpha \in V_P : \alpha \not \in \boldsymbol{h} \right\}$; \quad $V_P \leftarrow V_P \cup V_P^{\mathrm{new}} \setminus V_P^{\mathrm{out}}$; \quad $E_P \leftarrow  \{ E_P(\alpha) : \alpha \in V_P\}$ \label{ln:polyhedron-intersectionint-19}
	\STATE ${P} \leftarrow \mathtt{graph}(V_P,E_P)$  \label{ln:polyhedron-intersectionint-20}
	\STATE {$\mathbf{H} \gets \mathtt{insert}(\mathbf{H},V_P^\mathrm{new}); \quad \mathbf{H} \gets \mathtt{remove}(\mathbf{H},V^{\mathrm{out}}); \quad \mathbf{H} \gets \mathbf{H} \setminus \boldsymbol{h}$ }
	\label{ln:polyhedron-intersectionint-21}
	\RETURN $({P},\mathbf{H})$ \label{ln:polyhedron-intersectionint-22}
	\end{algorithmic}
	\label{alg:polyhedron-intersection}
	\end{algorithm}
 
\begin{figure}[http]
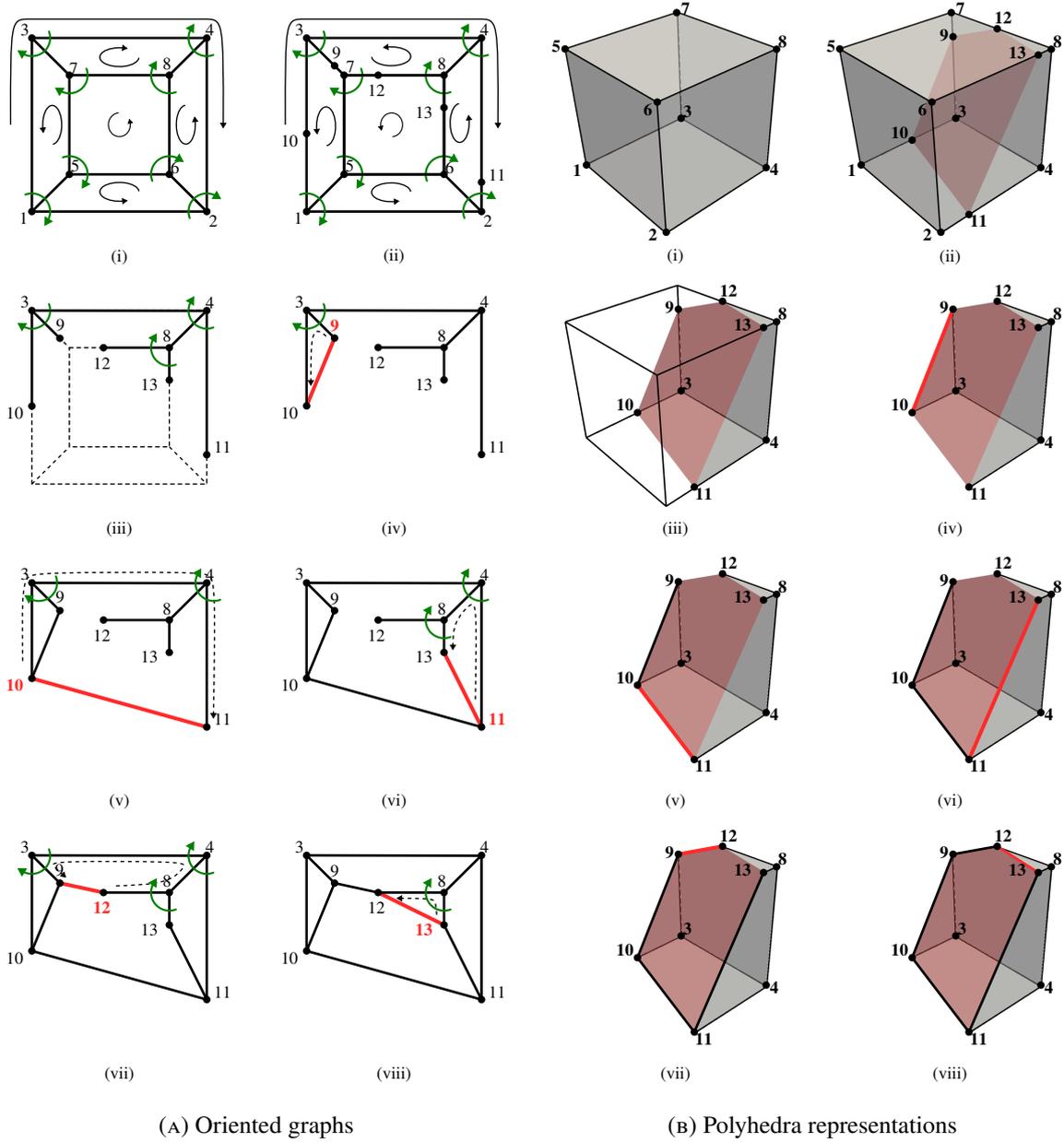

  \begin{subfigure}{0.49\textwidth}
    \centering
    \includefig[\tiny]{\textwidth}{alg2_graph}
    \caption{Oriented graphs}
  \end{subfigure}
  \begin{subfigure}{0.49\textwidth}
    \centering
    \includefig[\tiny]{\textwidth}{alg2_3d}
    \caption{Polyhedra representations}
  \end{subfigure}
  \caption{Illustration of \alg{alg:polyhedron-intersection} for an example in which we intersect a polyhedron, step (i), by a half-space plane. On the left-hand side figure (a), we show the polyhedron as a rotation system and how the intersection algorithm modifies this graph during different steps of the algorithm. On the right-hand side figure (b), we show the geometrical representation of the same steps. First, the polyhedron of step (i) is defined as a rotation system as described in line~\ref{ln:polyhedron-intersectionint-1} of the algorithm. Its edges are intersected by the half-space in step (ii) and the intersection points are computed (see line~\ref{ln:polyhedron-intersectionint-6}). Those new vertices are connected to the inside endpoints of the cut edges in step (iii) as indicated in line~\ref{ln:polyhedron-intersectionint-7}. In steps (v) to (viii), which correspond to the loop in line~\ref{ln:polyhedron-intersectionint-12}, the new vertices are connected to each other (see line~\ref{ln:polyhedron-intersectionint-17}) following an anti-clockwise path. The resulting polyhedron in step (viii) is represented with a new rotation system with all vertices inside the half-space.}
  \label{fig:alg2-simple-intersection}
\end{figure}
\subsection{Intersecting a polyhedron with a surface}\label{sec:intersection-P-S}

If the surface $S$ we want to intersect with the polyhedron $P$ is also convex, one can simply use Algorithm~\ref{alg:polyhedron-intersection} for all the half-spaces corresponding to the faces of $S$. However, $S$ is not convex for general geometries. In the final algorithm, we want to intersect an open oriented surface $\vec{\Gamma}$ (or its corresponding polyhedron representation $S \gets \mathtt{poly}(\vec{\Gamma})$) and a background mesh cell $K$. {We note that, for non-convex geometries, it is not possible to avoid the appearance of cells intersected with non-convex surfaces by using standard refinement strategies.}

In order to deal with general geometries, we perform a basic decomposition of the surface and the polyhedron into convex pieces~\cite{Chazelle1984}. 
First, we split the surface mesh $S$ into quasi-convex patches.  
In Algorithm~\ref{alg:reflex-walls} we compute a set of half-spaces that define such decomposition. 
{If $S$ is already a quasi-convex surface, the polytope is quasi-convex too, and we can proceed with the convex intersection as indicated above.}
{Otherwise, we identify cells intersected by non-convex surfaces by identifying the \emph{reflex edges} of the surfaces.} A reflex edge is the one that connects faces that are not quasi-convex. They are determined by a dihedral angle larger than $\pi$. Using the half-space representation, the reflex edges can be simply determined by the signed distance matrix block $\mathbf{H}_{S,S}$, i.e., the distance of $\mathtt{vert}(S)$ to the planes containing $\mathtt{faces}(S)$. In line~\ref{ln:reflex-walls-2} we iterate over all edges and extract the faces that share each edge in line~\ref{ln:reflex-walls-3}. An edge is reflex if the vertices of one face are exterior to the plane determined by the other face. We use this criterion to determine reflex edges in line ~\ref{ln:reflex-walls-4}. For quasi-coplanar faces, this definition can depend on the face being used to determine the plane. Besides, vertices of one face can lie on both sides of the plane, due to inexact arithmetic. In any case, these situations are not problematic when using the representations and algorithms discussed above. When two faces are quasi-aligned, considering the edge as reflex or not produces an error on the order of the machine precision. Besides, as discussed below, we consider a merging strategy to make quasi-aligned planes aligned in inexact arithmetic.

Standard methods to \emph{convexify} (i.e., split a non-convex polyhedron into convex parts) rely on vertical \emph{reflex walls}, i.e., vertical planes that contain the reflex edge. This algorithm is also denoted as \emph{vertical decomposition}. See~\cite{Hachenberger2008} for the application of this method to polyhedra. In Algorithm~\ref{alg:reflex-walls}, we do not consider a vertical wall. The definition of the plane without using information of the surface mesh is not a good choice in our case. We aim at reducing the intersections of the surface itself against these reflex walls and try to avoid quasi-aligned planes. Therefore, we consider the bisector of the two planes containing the faces sharing the reflex edge as the reflex wall (a \emph{bisection wall}). Thus, we compute this plane for each reflex edge in line~\ref{ln:reflex-walls-5} of Algorithm~\ref{alg:reflex-walls}. The result of this process for a surface $S$ (in polyhedral form) is represented with $\mathtt{walls}(S,\mathbf{H}_{S,S})$. 

\begin{algorithm}
  \caption{$\mathtt{walls}(S,\mathbf{H}_{S,S})$}
  \justifying
  \begin{algorithmic}[1]
    \STATE $R \gets \emptyset$
    \FOR{$ e \in \mathtt{edges}(S)$}  \label{ln:reflex-walls-2}
      \STATE $(T,U) \gets \mathtt{faces}(e)$ \label{ln:reflex-walls-3}
      \IF{$ \neg\left( \mathtt{vert}(U) \subset \mathbf{H}_{T,S} \, \land \, \mathtt{vert}(T) \subset \mathbf{H}_{U,S} \right) $} \label{ln:reflex-walls-4}
        \STATE $R \gets R \cup \mathtt{bisector}(T,U)$ \label{ln:reflex-walls-5} 
      \ENDIF
    \ENDFOR
    \RETURN $R$
  \end{algorithmic}
  \label{alg:reflex-walls}
\end{algorithm}

In Algorithm~\ref{alg:convex-decomposition} we decompose the surface $S$ and the polyhedron $P$ at hand into convex parts via the recursive splitting of these by the bisection walls of $S$, i.e., $\mathtt{walls}(S)$. In order to perform this algorithm, we need to compute first the signed matrix distance $\mathbf{H}_{SW,P}$. The first index $S$ stands for $\mathtt{faces}(S)$ planes while $W$ for $\mathtt{walls}(S)$ planes. The second index $P$ stands for $\mathtt{vert}(P)$ while $S$ stands for $\mathtt{vert}(S)$. We start the algorithm with $(\emptyset,(K,\mathbf{H}_{SW,K}),(S,\mathbf{H}_{SW,S}))$, where $K$ is a cell in the background mesh and $S$ the part of the whole surface mesh in touch with $K$. The recursivity in Algorithm~\ref{alg:convex-decomposition} is illustrated in Figure~\ref{fig:convex-decomposition}, where the decomposition of both surface $S$ and $K$ by the corresponding walls lead to a tree of pairs of convex surface and polyhedra components.}

Algorithm~\ref{alg:convex-decomposition} recursively intersects a polyhedron $P$ and surface $S$ against the walls and returns pairs of convex polytopes and surfaces after these intersections. If in the call to this recursive function there are still walls to be processed, we recursively \emph{convexify} $P$ and $S$ against each wall in lines~\ref{ln:convex-decomposition-7}-\ref{ln:convex-decomposition-8}. Note that we use closed half-space definitions for these intersections and that we convexify both sides after the intersection since both sides are of interest. We note that we can use either open or closed half-spaces in the intersections in lines~\ref{ln:convex-decomposition-7}-\ref{ln:convex-decomposition-8}.

If the function is invoked with no walls, we stop the process, since we have reached the leafs of the tree. The surface component $S$ that has been generated at this stage can still be disconnected, but what can be proved is that the connected components of $S$ are convex. If $S$ has disconnected components, we have to colour the surface into parts provide well-defined interiors of $P$. We do that in Algorithm~\ref{alg:colouring}. The reasoning behind this colouring is illustrated in Figure~\ref{fig:colouring}. We use these parts to colour $S$ in line~\ref{ln:convex-decomposition-3}. We return tuples of $P$ and each colour restriction of the surface $S$ in line~\ref{ln:convex-decomposition-4}. By construction, the interior of $P$ with respect to $S$ in each tuple is convex.     

\begin{algorithm}
	\caption{$\mathtt{convexify}(\boldsymbol{\mathcal{C}},
	({P},\mathbf{H}_{SW,P}),({S},\mathbf{H}_{SW,S}))$, }
	\begin{algorithmic}[1]
    \STATE $\mathbf{H}_{R,PS} \gets [\mathbf{H}_{R,P}, \mathbf{H}_{R,S}]$
	      \IF{$\mathbf{H}_{R,PS} = \emptyset$}
    \STATE $S \gets \mathtt{colouring}(S,\mathbf{H}_{S,S})$ \label{ln:convex-decomposition-3}
		\RETURN $\boldsymbol{\mathcal{C}} \leftarrow \boldsymbol{\mathcal{C}} \cup \left\{ ((P,\mathbf{H}_{T,P}),T) : {T} \in \mathtt{colours}(S) \right\}$ \label{ln:convex-decomposition-4}
	      \ELSE 
		\FOR{$\boldsymbol{h} \in \mathbf{H}_{R,PS}$}
		  \STATE $\boldsymbol{\mathcal{C}}^+ \gets \mathtt{convexify}(\boldsymbol{\mathcal{C}},(P,\mathbf{H}_{SW,P}) \cap \overline{\boldsymbol{h}_{{P}}}, (S,\mathbf{H}_{SW,S}) \cap \overline{\boldsymbol{h}_{{S}}})$
      \label{ln:convex-decomposition-7}
		  \STATE $\boldsymbol{\mathcal{C}}^- \gets \mathtt{convexify}(\boldsymbol{\mathcal{C}},(P,\mathbf{H}_{SW,P}) \cap \overline{-\boldsymbol{h}_{{P}}}, (S,\mathbf{H}_{SW,S}) \cap \overline{-\boldsymbol{h}_{{S}}})$
      \label{ln:convex-decomposition-8}
		  \STATE $\boldsymbol{\mathcal{C}} \leftarrow \boldsymbol{\mathcal{C}} \cup \boldsymbol{\mathcal{C}}^+ \cup \boldsymbol{\mathcal{C}}^-$
		\ENDFOR
	      \ENDIF
	\end{algorithmic}
	\label{alg:convex-decomposition}
\end{algorithm}
  
\begin{figure}[http]
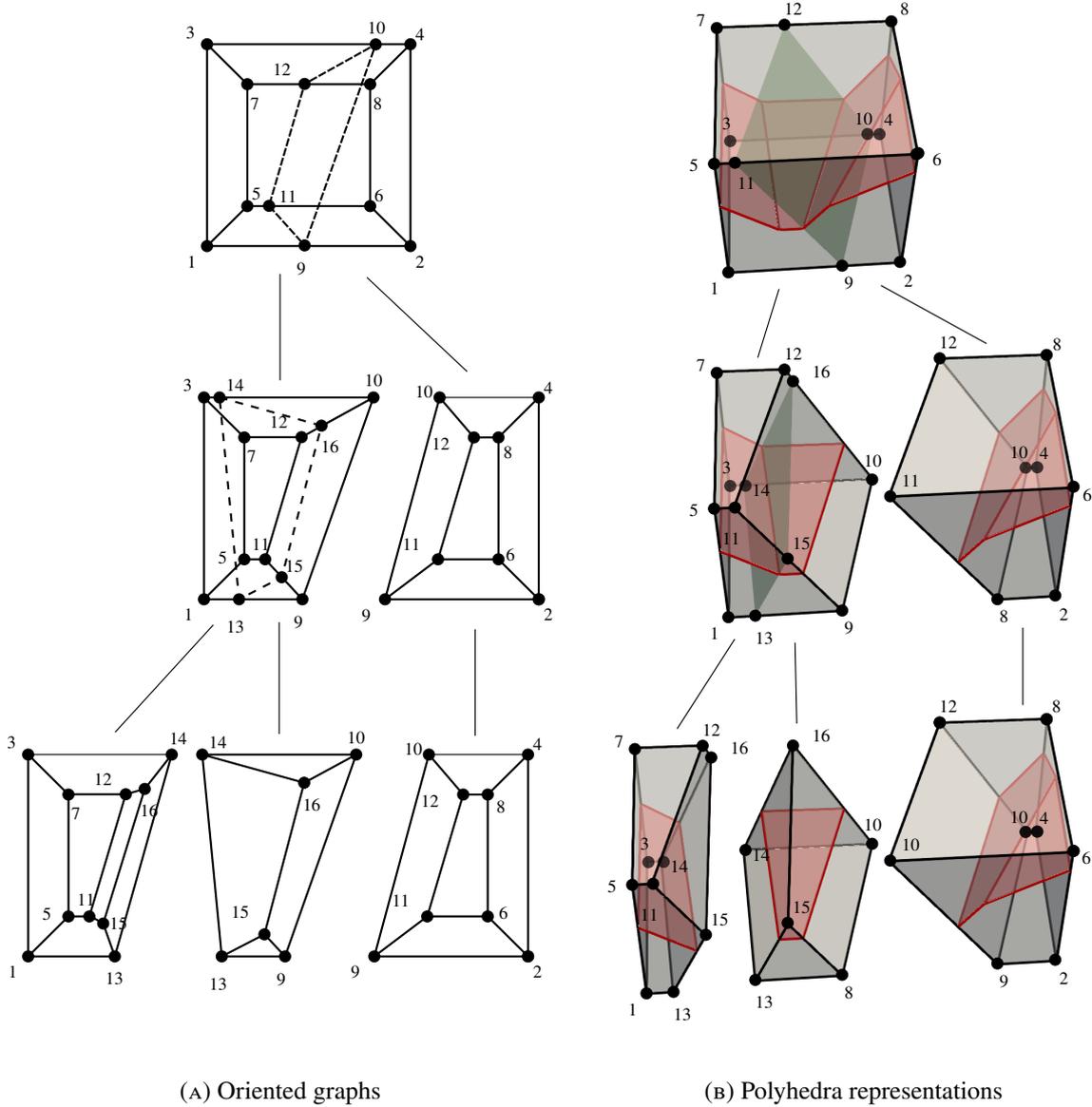

  \begin{subfigure}[b]{0.49\textwidth}
    \centering
    \includefig[\tiny]{\textwidth}{alg3_graph}
    \caption{Oriented graphs}
  \end{subfigure}
  \begin{subfigure}[b]{0.49\textwidth}
    \centering
    \includefig[\tiny]{\textwidth}{alg3_3d}
    \caption{Polyhedra representations}
  \end{subfigure}
  \caption{Example of the application of \alg{alg:convex-decomposition} to decomposes a cell polyhedron $K$ and a non-convex surface $S$ into convex parts. Both $K$ and $S$ are recursively split by the walls of $S$ using \alg{alg:polyhedron-intersection}. The clipping of $S$ is particularly simple since no extra vertex is introduced. Each row corresponds to a call of the algorithm. Recursion is introduced in line~\ref{ln:convex-decomposition-7}-\ref{ln:convex-decomposition-8} of the algorithm. The result is the leaves of the tree-like decomposition, which are processed in line~\ref{ln:convex-decomposition-3}-\ref{ln:convex-decomposition-4}. In each leaf, a piece of $K$ is associated with a convex piece of $S$.}
  \label{fig:convex-decomposition}
\end{figure}

\begin{figure}[http]
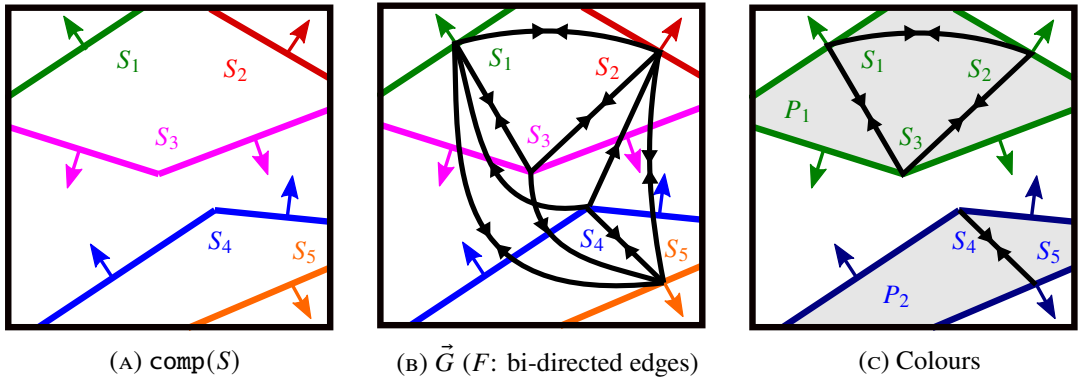

  \centering
  \begin{subfigure}{0.3\textwidth}
    \centering
    \includefig[\small]{.9\textwidth}{colouring-1}
    \caption{$\mathtt{comp}(S)$}
  \end{subfigure}
  \begin{subfigure}{0.3\textwidth}
    \centering
    \includefig[\small]{.9\textwidth}{colouring-2}
    \caption{$\vec{G}$ ($F$: bi-directed edges)}
  \end{subfigure}
  \begin{subfigure}{0.3\textwidth}
    \centering
    \includefig[\small]{.9\textwidth}{colouring-4}
    \caption{Colours}
  \end{subfigure}
  \caption{Illustration to explain Algorithm~\ref{alg:colouring} for a 2D example, which is called in line~\ref{ln:convex-decomposition-3} of Algorithm~\ref{alg:convex-decomposition}.  
  Given a surface $S$ with disconnected convex parts in (a), we build a directed graph $G$ (see (b)). $S_i$ is connected to $S_j$ if $S_i$ is inside $S_j$. Mutually connected components produce an undirected graph $F$. Let us discuss the first iteration of the while loop in line \ref{ln:colouring-while} for, e.g., $v$ being component $S_1$ of $S$. We find mutually connected components $NF$ to $S_1$ (including $S_1$) in line \ref{ln:colouring-nf-dg}. 
  $NF = \left\{ S_1, S_2, S_3, S_5 \right\}$ are the components that are in the interior of $S_1$ and $S_1$ is in their interior.  In the same line, we find the components $DG$ for which $S_1$ is outside. We get $DG = \left\{ S_4 \right\}$. We extract the mutually connected components to these ones, i.e., $\mathtt{adj}(F)(DG)$ and extract them from $NF$ in line \ref{ln:colouring-c-d}. We get $S_5$ and extract if from $NF$ to get the set $C = \left\{ S_1,S_2,S_3 \right\}$ that defines a convex polytope (have the same colour). We run the algorithm for the unprocessed components $D = \left\{ S_4, S_5 \right\}$ in the next iteration, which turn out to have the same colour. In this example, the algorithm returns two colours, namely $W = \{ (S_1,S_2,S_3), (S_4,S_5) \}$, which define two convex domains~$P_1$ and~$P_2$. We colour the graph $S$ with $W$ in line \ref{ln:colouring-colour}.  }
  \label{fig:colouring}
\end{figure}

\begin{algorithm}
	\caption{$\mathtt{colouring}(S,\mathbf{H}_{S,S})$}
	\begin{algorithmic}[1]
    \STATE $\vec{G} \gets \mathtt{graph}(\mathtt{comp}(S),(T,U) \to \mathtt{vert}(T) \subset \mathbf{H}_{U,S}) $ 
	  \STATE $F \gets \mathtt{graph}(\mathtt{comp}(S),(T,U) \to \mathtt{vert}(U) \subset \mathbf{H}_{T,S} \land \mathtt{vert}{(T)} \subset \mathbf{H}_{U,S} )$ 
    \label{ln:colouring-1}
    \STATE $V \gets \mathtt{vert}(\vec{G})$
    \STATE $W \gets \varnothing$
    \WHILE{$V \neq \emptyset$ \label{ln:colouring-while}}
       \STATE $v \gets V[1], \qquad NF \gets v \cup \mathtt{adj}(F)(v), \qquad DG \gets V \setminus (v \cup \mathtt{adj}(\vec{G})(v))$ \label{ln:colouring-nf-dg}
       \STATE $C \gets NF \setminus \mathtt{adj}(F)(DG), \qquad D \gets V \setminus C$ \label{ln:colouring-c-d}
       \STATE $W \gets W \cup \{C\}, \quad V \gets D$ \label{ln:colouring-7}
    \ENDWHILE 
    \STATE $S \gets \mathtt{colour}(S,W)$\label{ln:colouring-colour}
    \RETURN $S$
	\end{algorithmic}
	\label{alg:colouring}
\end{algorithm}

After Algorithm~\ref{alg:convex-decomposition}, we have a set of pairs of convex polyhedra and surfaces $(P,S)$. We can now use Algorithm~\ref{alg:polyhedron-intersection} for intersecting $P$ against all the half-spaces related to the faces of $S$. We do this in Algorithm~\ref{alg:intersection-P-S}. $\mathbf{H}$ must include the signed distance between vertices in $P$ and the half-spaces in $S$. The definition of open or close half-spaces depends on the definition of $S$ (as closed set, open set, or a mixed situation). At each leaf of the tree in Figure~\ref{fig:convex-decomposition}, Algorithm~\ref{alg:intersection-P-S} intersects $P$ as in Figure~\ref{fig:intersection-P-S}.

\begin{algorithm}
      \caption{$(P,\mathbf{H}) \cap {S} \to (P,\mathbf{H})$}
      \begin{algorithmic}[1]
        \FOR{$\boldsymbol{h} \in \mathbf{H}_{{S}}$} \label{ln:intersection-P-S-1}
            \STATE $(P,\mathbf{H}) \leftarrow (P,\mathbf{H}) \cap \boldsymbol{h}$ \label{ln:intersection-P-S-2}
            \ENDFOR
            \RETURN{$(P,\mathbf{H})$}
      \end{algorithmic}
      \label{alg:intersection-P-S}
\end{algorithm}

\begin{figure}[http]
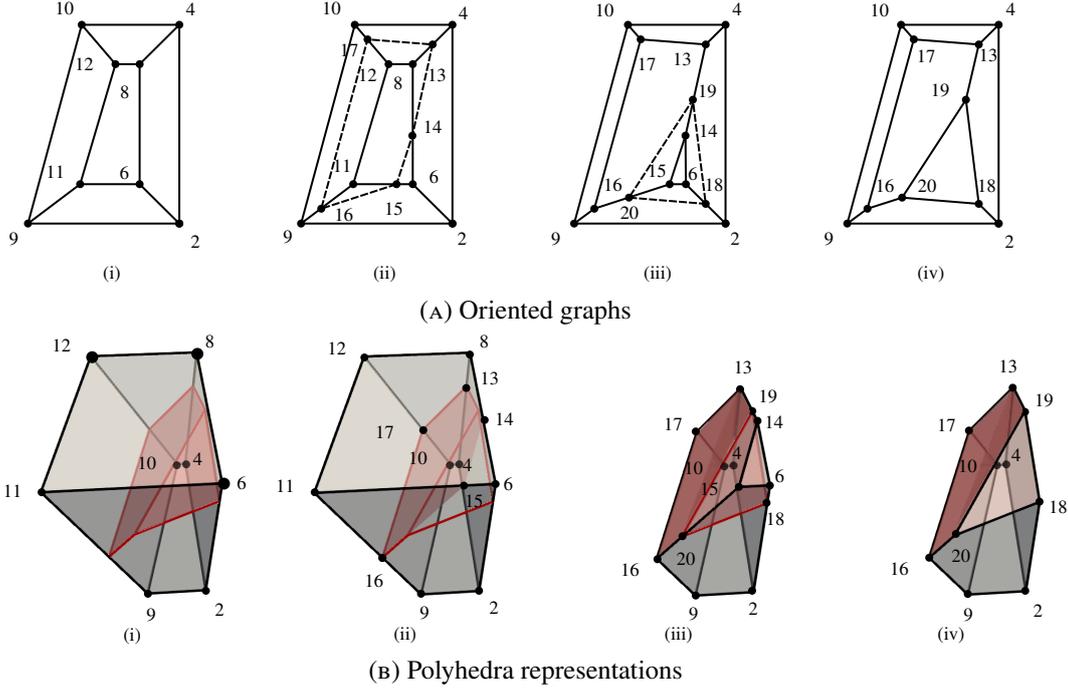

  \begin{subfigure}[b]{\textwidth}
    \centering
    \includefig[\tiny]{0.9\textwidth}{alg4_graph}
    \caption{Oriented graphs}
  \end{subfigure}
  \begin{subfigure}[b]{\textwidth}
    \centering
    \includefig[\tiny]{0.9\textwidth}{alg4_3d}
    \caption{Polyhedra representations}
  \end{subfigure}
  \caption{Illustration of \alg{alg:intersection-P-S}, which intersects  a convex volume $P$ by a convex surface $S$. The input of step (i) represents one of the leafs of \fig{fig:convex-decomposition}. 
  In steps (ii) and (iii) the volume polyhedra $P$ is intersected by the half-spaces determined by $\mathtt{faces}(S)$ as in the loop of line~\ref{ln:intersection-P-S-1}.
  The result in (iv) is the portion of $P$ inside $S$. If this process is repeated for each leaf of \fig{fig:convex-decomposition}, the result is the portion of $P$ in the interior defined by $S$, which is $P \cap \mathcal{B}$. }
  \label{fig:intersection-P-S}
\end{figure}

\subsection{Robust computation of signed distances}\label{sec:quasi-coplanarity}

The main problem when running the previous algorithms in inexact arithmetic is the computation of the signed distance matrices when multiple planes are quasi-aligned. For complex geometries, the number of surface mesh faces intersecting a background cell can still be large. There is a chance that some of these faces and the respective walls will be quasi-aligned. Even though this is not an issue in exact arithmetic, it can be problematic in inexact arithmetic. In this section, we provide mechanisms to enforce quasi-aligned half-spaces to be exactly aligned. Since we only make use of the signed distance matrix in our algorithms, the objective is to enforce the same entries for rows in $\mathbf{H}$ related to quasi-identical planes (with the same orientation) or times -1 for quasi-complimentary planes.

In a first step, we execute an algorithm $\mathtt{dist}(\Pi,V)$ that returns the the signed distance matrix $\mathbf{H}$ after computing the plane to vertex signed distances using parametric representations. In this method, the distances are \emph{snapped distances}. We define a \emph{snap} tolerance $\epsilon_{\mathrm{sn}}$ (e.g., 100 times the machine precision) and any distance within this tolerance is enforced to be zero. Geometrically, vertices extremely close to a plane are enforced to be on the plane in the half-space representation. It can happen that a node is snapped to multiple planes. We note that the snapping only affects the half-space representations; we do not perturb the vertices positions. 

In order to make the algorithm more robust, we additionally provide a mechanism to identify half-spaces that are quasi-identical or quasi-complimentary and to make them exactly aligned in the discrete representation. Algorithm~\ref{alg:align-surface} merges (aligns) discrete level-set representations of half-spaces $S$ (that represent surface faces) to the ones of a cell $K$ of the background mesh if they are quasi-aligned. The $S$ half-spaces, i.e., $\mathbf{H}_S$, can be perturbed in this process, but not the ones in $K$. Besides, the algorithm has been designed in such a way that the alignment of a half-space $\boldsymbol{h}^S$ against a half-space $\boldsymbol{h}^K$  of $K$ is consistent among all cells containing $\boldsymbol{h}^K$. 
In line~\ref{ln:align-surface-4} we check whether a $K$ half-space and an $S$ half-space are quasi-aligned. A surface half-space $\boldsymbol{h}^S$ is aligned with a cell half-space $\boldsymbol{h}^K$ if the absolute value of their distance to all the surface vertices in $S$  (in their discrete level-set representation) are below a given tolerance $\epsilon_{\mathrm{hs}}$. The absolute value is used to align not only two half-spaces that are quasi-coplanar but also the ones that are quasi-complementary. If the spaces are quasi-aligned, we run line~\ref{ln:align-surface-5}. Given two quasi-aligned planes $\boldsymbol{h}^i$ and $\boldsymbol{h}^j$, $\mathtt{sign}(\boldsymbol{h}^i,\boldsymbol{h}^j)$ returns +1 if they are quasi-coplanar and -1 if they are quasi-complementary. We note that this computation is numerically well-posed, e.g., comparing the sign of the distance to the furthest point in the discrete representation of the half-spaces. Finally, the $S$ half-space is replaced by the $K$ half-space times the sign, in order to keep consistency among cells.

\begin{algorithm}
  \caption{$\mathtt{align\_surface}(\mathbf{H}_{SK,KS}) \to \mathbf{H}_{S,KS}$ }
  \begin{algorithmic}[1]
    \FOR{$F_K \in \mathtt{faces}(K)$} 
      \FOR{$F_S \in \mathtt{faces}(S)$}
      \STATE $\boldsymbol{h}^i \gets \mathbf{H}_{F_K,S}: \ \boldsymbol{h}^j \gets \mathbf{H}_{F_S,S}$
        \IF{$\mathtt{min}( |\mathtt{max.}(\boldsymbol{h}^i-\boldsymbol{h}^j)| , |\mathtt{max.}(\boldsymbol{h}^i +\boldsymbol{h}^j)| )\leq \epsilon_{\mathrm{hs}}$
        \label{ln:align-surface-4}
        }         
          \STATE $\mathbf{H}_{F_S,KS} \gets \mathtt{sign}( \boldsymbol{h}^i, \boldsymbol{h}^j) \cdot \mathbf{H}_{F_K,KS}$ 
        \label{ln:align-surface-5}
        \ENDIF
      \ENDFOR
    \ENDFOR
    \RETURN $\mathbf{H}_{S,KS}$ 
  \end{algorithmic}
  \label{alg:align-surface}
\end{algorithm}

Once we have aligned the surface half-spaces with the cell half-spaces, modifying $\mathbf{H}_{S,KS}$, we are in position to run the cell-wise intersection algorithms. But we still need to check whether wall and surface planes are quasi-aligned. In Algorithm~\ref{alg:align-planes} we provide an algorithm that aligns planes given a general signed distance matrix $\mathbf{H}$. First, the algorithm creates in line~\ref{ln:align-planes-1} a graph of half-spaces in which two half-spaces are connected if they are quasi-aligned, using the same condition as in Algorithm~\ref{alg:align-surface} but for all vertices in the discrete representation of the half-spaces. We extract the components of these graphs in line~\ref{ln:align-planes-2}. Half-spaces in a component are considered to be all quasi-aligned and we enforce all their distances to be the same (times -1 for quasi-complementary half-spaces). The value of the distances that we have used is provided in Algorithm~\ref{alg:equal-rows}. In this algorithm, vertices that belong to one of the half-spaces in a component belong to the half-space after alignment, i.e., the distance is 0 (line~\ref{ln:equal-rows-7}). For other vertices, we just pick the signed distance from one of the half-spaces in line~\ref{ln:equal-rows-8} (with the right sign, computed in line~\ref{ln:equal-rows-4}). In any case, as soon as the merge is performed, other reasonable choices could also be considered without affecting the robustness of the overall algorithm. After this process, all quasi-aligned components are exactly aligned.

\begin{algorithm}
  \caption{$\mathtt{align\_planes}(\mathbf{H}) \to \mathbf{H}$ }
      \justifying
  \begin{algorithmic}[1]
    \STATE $G \gets \mathtt{graph}(\mathbf{H}, (\boldsymbol{h}^i,\boldsymbol{h}^j) \to \mathtt{min}( |\mathtt{max.}(\boldsymbol{h}^i-\boldsymbol{h}^j)| , |\mathtt{max.}(\boldsymbol{h}^i +\boldsymbol{h}^j)| )\leq \epsilon_{\mathrm{hs}}$)
    \label{ln:align-planes-1}
     \STATE $C \gets \mathtt{comp}(G)$
     \label{ln:align-planes-2}
        \FOR{$T \in C$} 
        \STATE $\mathtt{merge}(\mathbf{H}_{T,*})$
        \ENDFOR 
    \RETURN $\mathbf{H}$ 
  \end{algorithmic}
  \label{alg:align-planes}
\end{algorithm}

\begin{algorithm}
  \caption{$\mathtt{merge}(\mathbf{H}) \to \mathbf{H}$ }
      \justifying
  \begin{algorithmic}[1]
\STATE $\boldsymbol{h}^0 \gets \mathbf{H}{[1,:]}$
\STATE $\boldsymbol{s} \gets \mathtt{zeros}(\mathtt{dims}(\mathbf{H}){[1]})$
\FOR{$(i,\boldsymbol{h}) \in \mathtt{enum}(\mathtt{rows}(\mathbf{H}))$}
\STATE $\boldsymbol{s}[i] \gets \mathtt{sign}(\boldsymbol{h},\boldsymbol{h}^0)$
   \label{ln:equal-rows-4}
\ENDFOR
\FOR{$\boldsymbol{c} \in \mathtt{columns}(\mathbf{H})$}
   \STATE $(\mathtt{min}(\mathtt{abs.}(\boldsymbol{c})) = 0) \ ? \ d \gets 0 : d \gets \boldsymbol{c}[1]$
   \label{ln:equal-rows-7}
   \STATE $\boldsymbol{c} \gets d \cdot \boldsymbol{s}$
   \label{ln:equal-rows-8}
    \ENDFOR
    \RETURN $\mathbf{\textbf{H}}$ 
  \end{algorithmic}
  \label{alg:equal-rows}
\end{algorithm}

\subsection{Global intersection algorithm}\label{sec:global}

We are in position to define Algorithm~\ref{alg:global}, the global algorithm we propose to intersects a background mesh $\mathcal{T}$ and a boundary mesh $\mathcal{B}$. The results is a partition of each cell in both meshes into sub-cells, denoted with $\mathcal{T}^{\mathrm{cut}}, \, \mathcal{B}^{\mathrm{cut}}$. Figure~\ref{fig:global} illustrates all the steps being performed in this algorithm to intersect a cell $K\in\mathcal{T}$ with the boundary mesh $\mathcal{B}$.

First, we perform a background cell-wise intersection (see line~\ref{ln:global-3}). In general, the surface mesh can have a large number of cells but a background cell usually intersects a very small portion of these surface cells. For computational efficiency and robustness of the algorithm, it is essential to reduce the polyhedron clipping to the portion of the surface $\mathcal{B}$ that can be in touch with the cell. This step is denoted with $\mathtt{restrict}$ in line~\ref{ln:global-3}. It makes use of cheap geometrical predicates, since the result does not need to be precise; false positives do not pose any problem. In fact, in order to capture cells that are quasi-aligned to faces in the background cell $K$, we need to enlarge $K$  at least a distance equal to $\epsilon_{\mathrm{hs}}$. Since these predicates are quite standard in computational geometry and can be found in computational geometry libraries like CGAL \cite{cgal:eb-21b}, they are not included here for the sake of conciseness.

After the restriction, we transform the portion of the surface mesh into a polyhedron in line~\ref{ln:global-4}, using Algorithm~\ref{alg:surface-to-polyhedron}. In line~\ref{ln:global-5}, we compute the signed distance matrix between the vertices in $\mathtt{vert}(K) \cup \mathtt{vert}(S)$ and the planes in $\mathtt{faces}(K) \cup \mathtt{faces}(S) \cup \mathtt{walls}(S)$ using standard algorithms.

The signed (snapped) distance matrix for all faces of $K$ and $S$, and walls of $S$, and vertices in $K$ and $S$ is computed in line~\ref{ln:global-5}. Next, we align the surface half-spaces to the ones of the cell boundaries in a consistent way using Algorithm~\ref{alg:align-surface} in line~\ref{ln:global-6}.

Given the rectangular cell $K = [x^{-},x^+] \times [y^{-},y^+] \times [z^{-},z^+]$, we define $K_{{\circ\bullet}} = (x^{-},x^+] \times (y^{-},y^+] \times (z^{-},z^+]$. In order to perform the surface mesh cell-wise intersection, we precisely use  $K_{{\circ\bullet}}$, not $K$, in line~\ref{ln:global-7} using Algorithm~\ref{alg:intersection-P-S}.  Otherwise, faces that lie on background cell boundaries would be processed twice. The result of this surface-cell intersection (line~\ref{ln:global-8}) for all cells returns a refinement of $\mathcal{B}$, denoted with $\mathcal{B}^{\mathrm{cut}}$.  Such intersection is illustrated in Figure~\ref{fig:intersect-S-K}.\footnote{We note that one could also extract the surface mesh as the surface of the clipped polytopes obtained after intersecting again half-spaces in $\mathtt{faces}(S)$ in line~\ref{ln:global-12}. This is the reason why we use $K_{\bullet\circ}$ in line~\ref{ln:global-7} (to process surface faces aligned with background cells faces only once), closed spaces in lines~\ref{ln:convex-decomposition-7}-\ref{ln:convex-decomposition-8} of Algorithm~\ref{alg:convex-decomposition} (not to lose any surface face after splitting with wall half-spaces) and intersection against open half-spaces related to the surfaces in line~\ref{ln:global-11} (to discard zero volume components after this decomposition and count surface faces on walls only once). A simple example in which we can encounter this situation is illustrated in Figure~\ref{fig:open-close}. In any case, these choices of open/closed half-spaces are not required when extracting the surface mesh as in line~\ref{ln:global-8}. The connection between the interior partition $\mathcal{T}^{\mathrm{cut}}$ and $\mathcal{B}^{\mathrm{cut}}$ is not important for the unfitted scheme being used later on because Dirichlet boundary conditions are weakly imposed. In any case, it could be useful for other embedded methods that would make use of a strong imposition of Dirichlet data.} We can readily use $\mathcal{B}_D$ and $\mathcal{B}_N$ instead, to compute $\mathcal{B}^{\mathrm{cut}}_D$ and $\mathcal{B}^{\mathrm{cut}}_N$.

\begin{figure}[http]
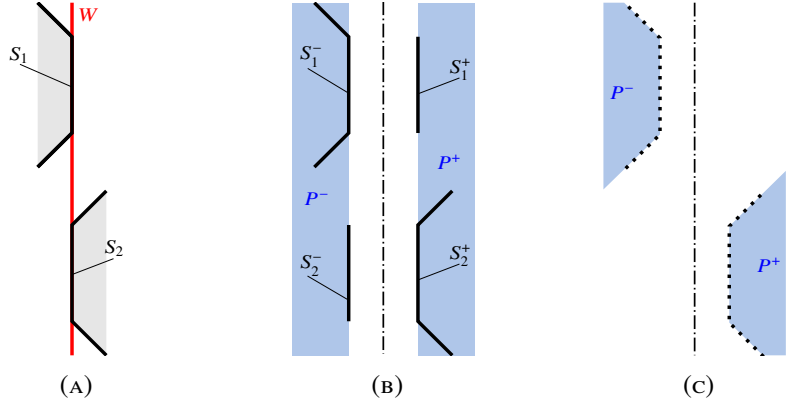

  \centering
  \begin{subfigure}{0.25\textwidth}
    \centering
    \includefig{0.6\textwidth}{open_close_1}
    \caption{}
  \end{subfigure}
  \begin{subfigure}{0.25\textwidth}
    \centering
    \includefig{0.6\textwidth}{open_close_2}
    \caption{}
  \end{subfigure}
  \begin{subfigure}{0.25\textwidth}
    \centering
    \includefig{0.6\textwidth}{open_close_3}
    \caption{}
  \end{subfigure}
  \caption{Simple 2D example to justify the choice of open and closed half-spaces in Algorithm~\ref{alg:convex-decomposition} and \ref{alg:global}. In (a) $S_1$ and $S_2$, surface faces are aligned with the wall $W$. When decomposing by $W$ with closed half-spaces, i.e., $-\overline{\boldsymbol{h}_W}$ and $\overline{\boldsymbol{h}_W}$ in lines~\ref{ln:convex-decomposition-7}-\ref{ln:convex-decomposition-8} of Algorithm~\ref{alg:convex-decomposition}, $S_1$ and $S_2$ are repeated in both sides as it is shown in (b). However, when intersecting $P\cap S$, since $\boldsymbol{h}_S$ is open in line~\ref{ln:global-11} of Algorithm~\ref{alg:global}, $S_1^+$ (resp., $S_2^-$) does not belong to the open half-space $\boldsymbol{h}_{S_1}$ (resp., $\boldsymbol{h}_{S_2}$) and thus eliminated after the intersection. The resulting polyhedra after clipping are shown in (c). This specific choice is required in the case in which one wants to extract the boundary surface from the clipped polytopes, i.e., $\mathcal{B}^{\mathrm{cut}} \gets \partial \mathcal{T}^{\mathrm{cut}}$.}
  \label{fig:open-close}
\end{figure}

In this intersection, since the half-spaces related to $\mathtt{faces}(K)$ are processed, they are eliminated from $\mathbf{H}$ (see Algorithm~\ref{alg:polyhedron-intersection}). Finally, we merge quasi-aligned surface and wall half-spaces using Algorithm~\ref{alg:align-planes}. With the resulting signed distance matrix, we run the convex decomposition in Algorithm~\ref{alg:convex-decomposition} in line~\ref{ln:global-10}, starting with the polyhedron $K$ and surface $S$. {We note that this line is doing nothing if there are no walls, i.e., if the polytope is already quasi-convex.} The resulting convex polyhedron-surface components are intersected using Algorithm~\ref{alg:polyhedron-intersection} and added to the sub-mesh for $K$ that represent its interior part in line~\ref{ln:global-11}. 

\begin{algorithm}
	\caption{$\T \cap \B \to \mathcal{T}^{\mathrm{cut}},\B^{\mathrm{cut}} $ }
	\justifying
	\begin{algorithmic}[1]
	  \STATE $\Tcut \gets \emptyset; \quad \Bcut \gets \emptyset$,
    \FOR{$K \in \T$} \label{ln:global-2}
	    \STATE $B \leftarrow \mathtt{restrict}(\mathcal{B},K)$ \label{ln:global-3} 
	    \STATE $S \gets \mathtt{pol}(B)$ \label{ln:global-4}
	    \STATE $\mathbf{H}_{KSW,KS} \gets \mathtt{dist}( [ \mathtt{faces}(K), \mathtt{faces}(S), \mathtt{walls}(S) ] ,[ \mathtt{vert}(K) , \mathtt{vert}(S) ])$ \label{ln:global-5} 
	    \STATE $\mathbf{H}_{S,KS} \gets \mathtt{align\_surface}(\mathbf{H}_{S,KS},\mathbf{H}_{K,KS})$ \label{ln:global-6}
	    \STATE $(S,\mathbf{H}_{S,KS}) \gets (S,\mathbf{H}_{KS,KS}) \cap {K}_{\mathrm{\circ\bullet}}$ \label{ln:global-7}
	    \STATE $\Bcut \gets \Bcut \cup S$ \label{ln:global-8}
	    \STATE $\mathbf{H}_{SW,KS} \gets \mathtt{align\_planes}(\mathbf{H}_{SW,KS})$ \label{ln:global-9}
	    \STATE $\C \leftarrow \mathtt{convexify}(\emptyset,(K,\mathbf{H}_{SW,K}),(S,\mathbf{H}_{SW,S}))$ \label{ln:global-10}
	    \STATE $\mathcal{T}_K \gets \left\{ P \cap S : ((P,\mathbf{H}_{S,P}),S) \in \C \right\}; \quad \Tcut \gets \Tcut \cup \mathcal{T}_K$ \label{ln:global-11}
	  \ENDFOR\label{ln:global-12}
	  \RETURN{$\mathcal{T}^{\mathrm{cut}}, \, \mathcal{B}^{\mathrm{cut}}$} \label{ln:global-13}
	\end{algorithmic}
	\label{alg:global}
\end{algorithm}

\begin{figure}[http]
  \begin{subfigure}[b]{0.49\textwidth}
    \centering
    \includefig[\tiny]{\textwidth}{intersect_S_K_graph}
    \caption{Oriented graphs}
  \end{subfigure}
  \begin{subfigure}[b]{0.49\textwidth}
    \centering
    \includefig[\tiny]{\textwidth}{intersect_S_K_3d}
    \caption{Polyhedra representations}
  \end{subfigure}
  \caption{Illustration of \alg{alg:intersection-P-S}, i.e., $(S,\mathbf{H}) \cap K $, when calling line~\ref{ln:global-7} of \alg{alg:global}. We consider $S$ and $K$ from step (i) in \fig{fig:global}.
  In this process, step (i) to (v), $S$ is intersected by each half-space $\boldsymbol{h}_i \in \mathbf{H}_{S,K}$ related to $\mathtt{faces}(K)$ using \alg{alg:polyhedron-intersection}, which is called in the loop of line~\ref{ln:intersection-P-S-1} of \alg{alg:intersection-P-S}. Note that $\boldsymbol{h}_3$ (bottom plane) is excluded from the figure because the intersection is meaningless. The result, step (vi), is a new surface $S$ inside $K$, which is introduced in step (ii) of \fig{fig:global}.}
  \label{fig:intersect-S-K}
\end{figure}

\begin{figure}[http]
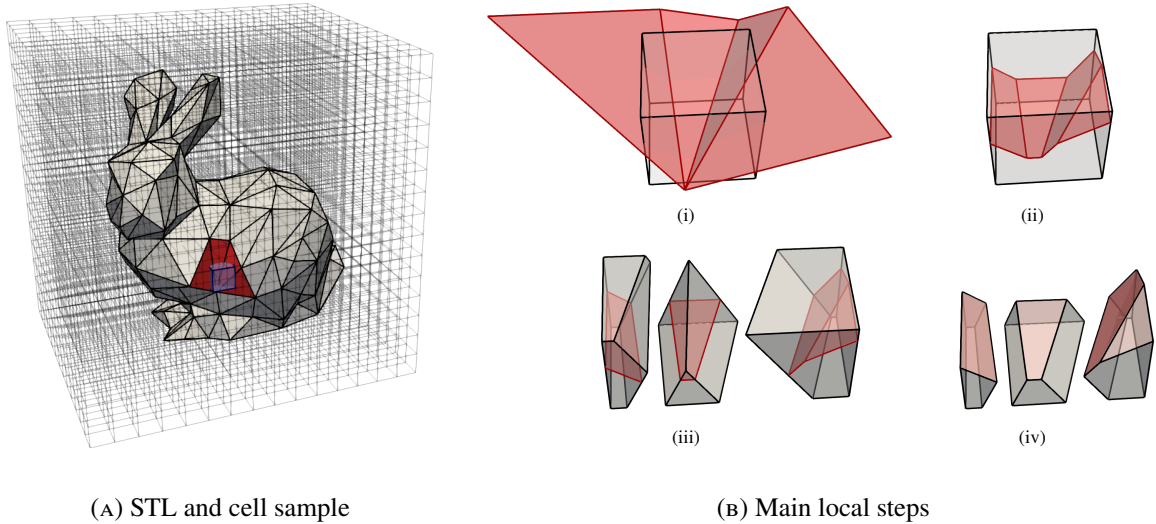

  \begin{subfigure}[b]{0.39\textwidth}
    \centering
    \includefig[\tiny]{\textwidth}{alg5_a}
    \caption{STL and cell sample}
  \end{subfigure}
  \begin{subfigure}[b]{0.59\textwidth}
    \centering
    \includefig[\tiny]{\textwidth}{alg5_b}
    \caption{Main local steps}
  \end{subfigure}
  \caption{Illustration of \alg{alg:global}, which, given an \ac{stl} $\mathcal{B}$ and a background mesh $\mathcal{T}$ (see (a)) intersects each background cell $K \in \mathcal{T}$ with $\mathcal{B}$. The steps described in (b) represent the loop in line~\ref{ln:global-2}. First,  $\mathcal{B}$ is restricted to the faces touching $K$ in (i) (line~\ref{ln:global-3}) and defined as a polyhedron $S$ (line~\ref{ln:global-4}). Next,  $S$ is intersected by the half-spaces bounding $K$ performed in step (ii) (see line~\ref{ln:global-7}). As the surface $S$ may be non-convex, it is decomposed into convex parts in step (iii), as described in line~\ref{ln:global-10}, together with $K$. Finally, in line~\ref{ln:global-11}, each convex component $P$ of $K$ is intersected by the corresponding part of $S$. The result in step (iv) is $K \cap \Omega$ described as the union of convex polyhedra, represented as in Definition~\ref{def:pol}.
  Steps (ii), (iii) and (iv) are further detailed in Figure~\ref{fig:intersect-S-K}, Figure~\ref{fig:convex-decomposition} and Figure~\ref{fig:intersection-P-S} respectively.} 
  \label{fig:global}
\end{figure}

Even though we have presented the algorithm the interior component only, i..e, $K \cup \Omega$, it is computationally efficient to compute the convex decomposition of both interior and exterior at the same time when the latter is needed, e.g., in interface problems. We also note that the definition of interior, exterior and boundary vertices is a straightforward side-result of the algorithm. The interior (resp., exterior) is determined in line~\ref{ln:global-11} for cut cells and propagated globally to other interior (resp., exterior) cells.

Some algorithms, e.g., the numerical integration, could require a simplex decomposition of the polyhedra that define the interior or exterior. A convex decomposition of a convex polytope is straightforward and can be computed symbolically (see~\cite{Powell2015} for details). In any case, this step is optional. Even for numerical computations, one can use quadrature rules for general polytopes that do not require this step~\cite{Chin2020}.

\section{Numerical experiments}\label{sec:num-exp}

\subsection{Objectives}

In the numerical examples below, we analyse the algorithmic and computational performance of the intersection algorithm proposed in this paper. In particular, we study the accuracy of the intersection method, its robustness,  the scaling of CPU times with respect to the number of cells in the background mesh and the faces of the \ac{stl} and its usage in unfitted \ac{fe} simulations. We consider three different numerical experiments. In the first one (\sect{sec:thingi10k}), we run the intersection algorithm in a large set of geometries taken form the Thingi10K \cite{Zhou2016} collection of \ac{stl} models. We apply the method to all models in this data-set that fulfil the input requirements of the intersection algorithm in order to evaluate its ability to deal  with complex and arbitrary inputs. In the second experiment (\sect{sec:robustness}), we analyse the robustness of the method with respect to perturbations in the background mesh, either with translations or rotations. And finally (\sect{sec:fe-test}), we apply the proposed intersection method to generate  integration cells in an unfitted \ac{fe} method to analyse the influence of the cutting algorithm in the quality of the \ac{fe} solution. We have performed the simplex decomposition step in all experiments.

\subsection{Experimental setup} 
The numerical experiments have been performed on TITANI, a medium size cluster at the Universitat Politècnica de Catalunya (Barcelona, Spain) and on Gadi, a high-end supercomputer at the NCI (Australia) with 3024 nodes, each one powered by a 2 x 24 core Intel Xeon Platinum 8274 (Cascade Lake) at 3.2 GHz and 192GB RAM. The timing experiments have been performed on Gadi exclusively, whereas TITANI has been considered for non-performance critical runs. In order to reduce the influence of external factors on the CPU timings, each time measure reported in the experiments is computed as the minimum of 5 runs in the same Julia session, i.e., one run for Julia JIT compilation and four runs to measure run-time performance. The intersection algorithms have been implemented using the Julia programming language \cite{Julia-2017} and are freely available in the STLCutters.jl package \cite{Martorell_STLCutters_2021}. The unfitted \ac{fe} computations have been performed using the Julia \ac{fe} library Gridap.jl \cite{Badia2020b} version 0.16.3 and the extension package for unfitted methods GridapEmbedded.jl~\cite{GridapEmbedded-jl} version 0.7. In order to parse the \ac{stl} files, we have used the MeshIO.jl \cite{MeshIO-jl} Julia package version 0.4.

\subsection{Batch processing the \ac{stl} models of the Thingi10K data-set}
\label{sec:thingi10k}

We start the numerical experiments by processing a large number of real-world \ac{stl} models to show the capacity of the proposed intersection algorithm to deal with complex and arbitrary data automatically. To this end, we consider the Thingi10K \cite{Zhou2016} data-base, which contains ten thousand 3D \ac{stl} models, from simple to very complex, used mainly for real-world 3D printing purposes. Our goal is to show that our intersection algorithms are able to handle these geometries automatically and directly without any manual pre-process as a demonstration of the robustness and generality of the proposed algorithm. 
 
 Among all models within the Thingi10k set, we process the ones that fulfil the requirements of our method. In particular, we need closed surfaces that define a volume. Not all geometries in the database fulfil this condition and, thus, we extract valid geometries by considering the ones tagged as \emph{is closed} and \emph{is manifold}. E.g., one can recover these geometries by typing ``is closed, is manifold'' in the search field of the Thingi10k web page. This results in a subset of 4963 models. Among them, we have found 211 cases that could not be processed either due to broken download links or corrupt STL files (i.e., the parser was not able to read the model into memory) and 20 cases that are not a manifold up to machine precision. By discarding these pathological cases, we recovered the 4732 geometries that  have been processed in this test. As an example, \fig{fig:mat-geo} shows  some of the processed \ac{stl} models, which illustrates the diversity of cases analysed in this experiment.

\begin{figure}[http]
  \includegraphics[width=0.9\textwidth]{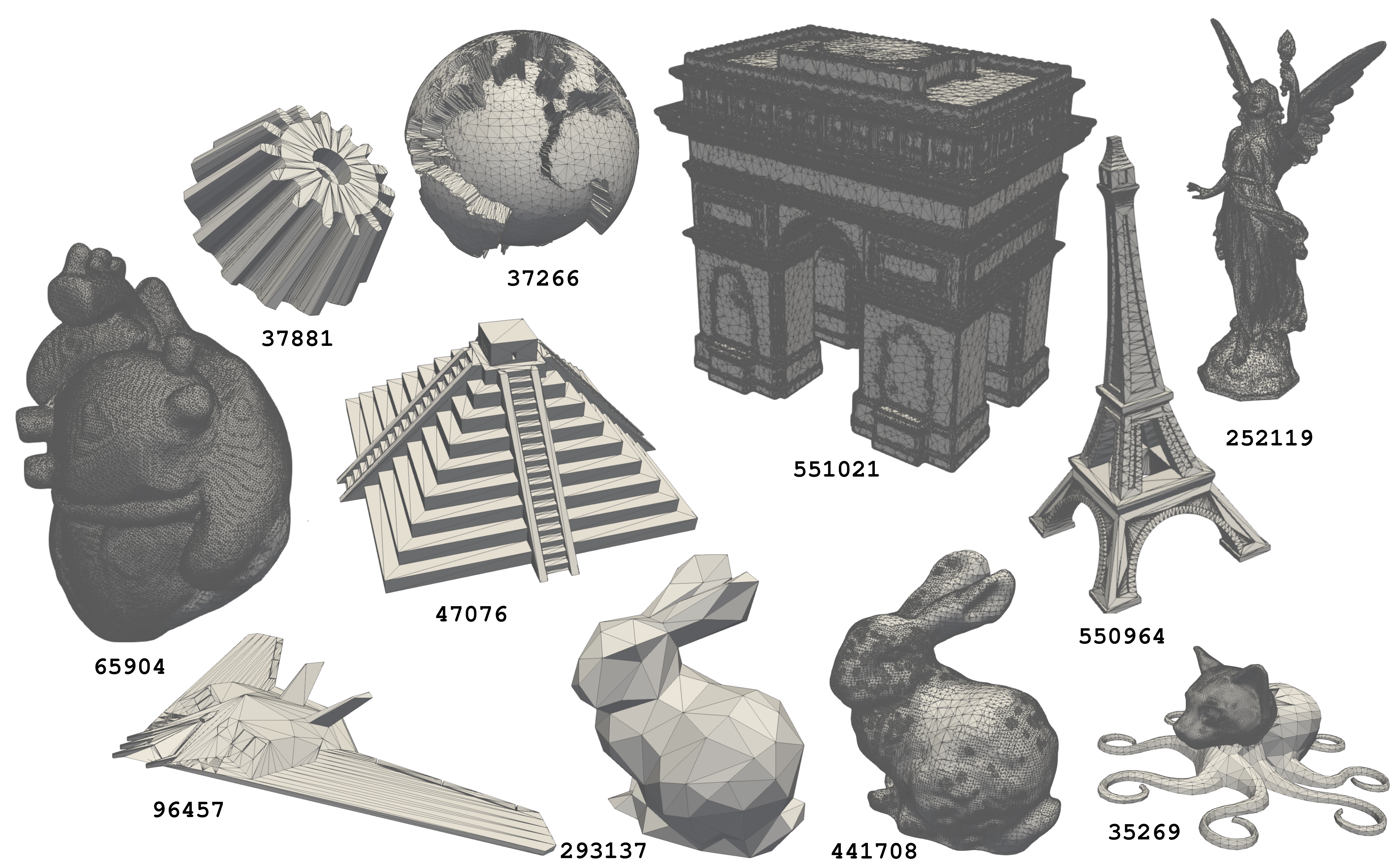}
  \caption{Selection of 11 \ac{stl} models from the Thingi10K database processed in the numerical examples. They are displayed with their corresponding model id provided by the Thingi10K database.}
  \label{fig:mat-geo}
\end{figure}

Each of the considered models is processed automatically as follows. First, we parse the downloaded \ac{stl} file and compute its bounding box extreme points $\mathbf{p}^{\rm min}_{\rm stl}$ and $\mathbf{p}^{\rm max}_{\rm stl}$. Then, we generate a 3D background Cartesian mesh, which covers a box approximately 40\% larger in each direction than the bounding box of the \ac{stl}. We generate the Cartesian mesh with at least $n^{\rm max}=100$ cells in the largest axis and $n^{\rm min}=10$ in the shortest. The element size $h$ and the bounding box points of the background Cartesian mesh, $\mathbf{p}^{\rm min}_{\rm msh}$ and $\mathbf{p}^{\rm max}_{\rm msh}$, are respectively computed  as
\begin{equation}
h = 1.4\min\left\{  \max\left( \dfrac{\mathbf{p}^{\rm max}_{\rm stl} - \mathbf{p}^{\rm min}_{\rm stl}}{n^{\rm max}} \right), \min\left( \dfrac{\mathbf{p}^{\rm max}_{\rm stl} - \mathbf{p}^{\rm min}_{\rm stl}}{n^{\rm min}} \right)     \right\},
\end{equation}
and
\begin{equation}
\mathbf{p}^{\rm min}_{\rm msh} \doteq \mathbf{p}^{\rm min}_{\rm stl} -0.2\left(\mathbf{p}^{\rm max}_{\rm stl} - \mathbf{p}^{\rm min}_{\rm stl} \right),\quad
  \mathbf{p}^{\rm max}_{\rm msh} \doteq  \mathbf{p}^{\rm min}_{\rm msh} + \left\lceil \dfrac{ 1.4 \left( \mathbf{p}^{\rm max}_{\rm stl} - \mathbf{p}^{\rm min}_{\rm stl} \right)}{h} \right\rceil h.
\end{equation}

\begin{figure}[http]
  \centering
  \begin{subfigure}{0.2\textwidth}
    \includegraphics[height=0.23\textheight]{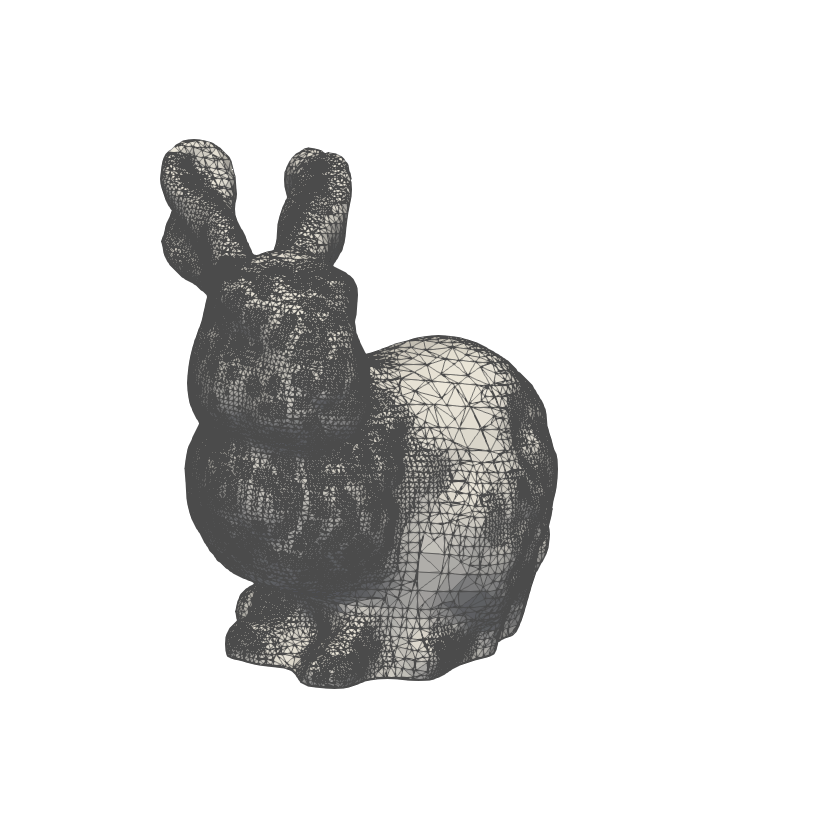}
    \caption{}
  \end{subfigure}
  \begin{subfigure}{0.5\textwidth}
    \includegraphics[height=0.23\textheight]{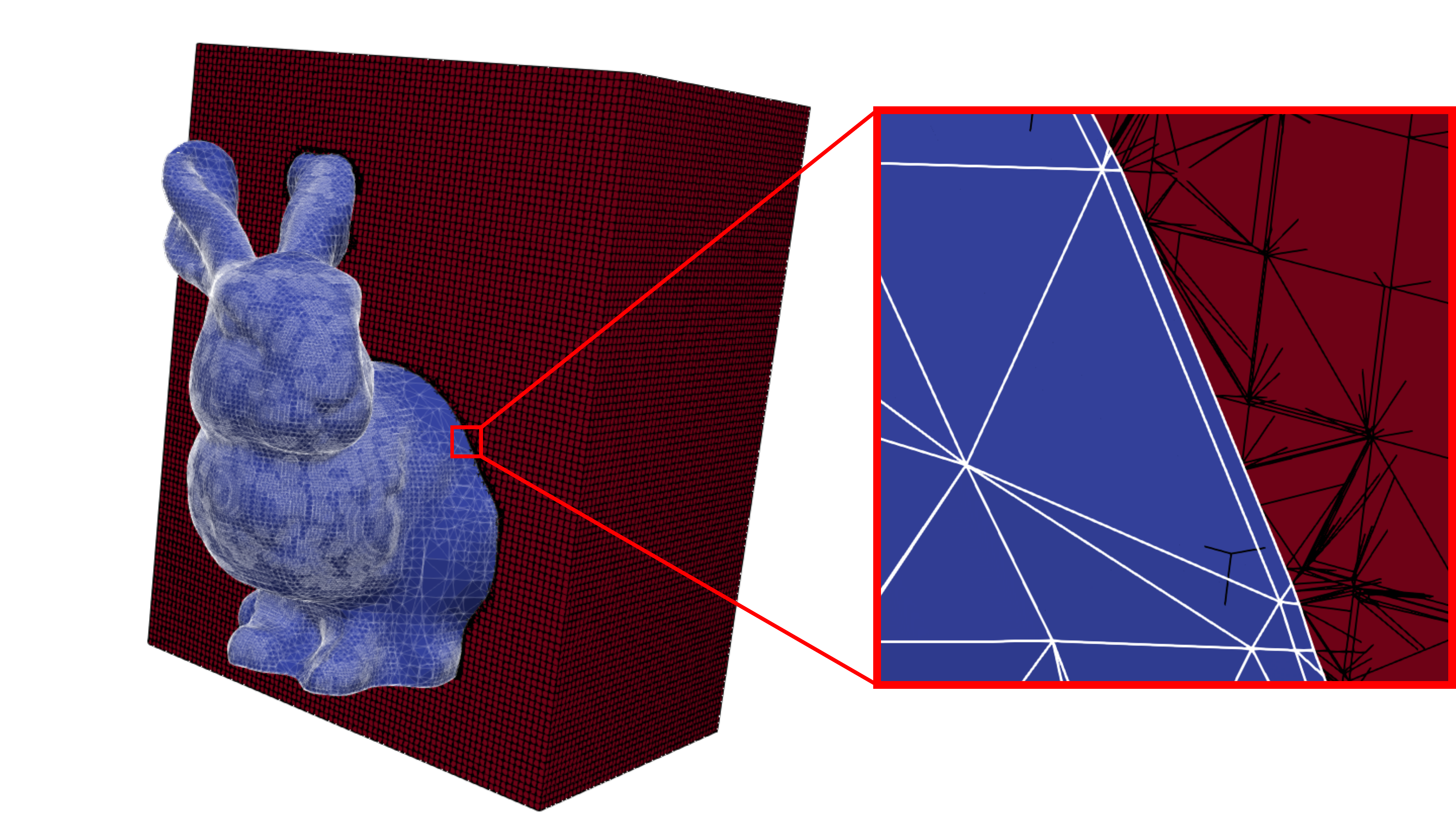}
    \caption{}
  \end{subfigure}
  \caption{Generated volume sub-triangulation for the \ac{stl} model with id 441708: (a) shows the original \ac{stl} geometry, while (b) shows a clipped portion of the volume sub-triangulation (red cells) with a detail of the cut cells near the \ac{stl} faces defining the boundary (blue faces).}
  \label{fig:bunny-stl+zoom}
\end{figure}

In a next step, the background mesh is intersected with the \ac{stl} surface mesh using \alg{alg:global}. See, e.g., in \fig{fig:bunny-stl+zoom} a detail of the resulting volume sub-triangulation for one of the considered geometries. The intersection algorithm is applied with snap tolerance $\esn = \ell^\mathrm{max}_\B 10^2\mathrm{eps}$ and quasi co-planar tolerance $\ehs = \ell^\mathrm{max}_\B 10^3\mathrm{eps}$,  being $\ell^\mathrm{max}_\B$ the length of the largest axis of the \ac{stl} bounding box and $\mathrm{eps}$ the machine precision associated with 64-bit floating point numbers. The final step is to compute some indicators of the quality of the generated sub-triangulations.  On the one hand, we measure $\Gamma^{\mathrm {st}}$, the area of the boundary sub-triangulation and compare it with $\Gamma^\mathrm{STL}$, the area of the original \ac{stl} mesh. From these values, we compute the relative surface error $\epsilon_\Gamma = | \Gamma^{\mathrm{STL}} - \Gamma^{\mathrm{st}} | / \Gamma^{\mathrm{STL}}$. As the input geometries are represented by surfaces, the original interior volume is unknown. Thus, we quantify the volume error by comparing the inside and outside volumes of the bulk sub-triangulation, $V^\mathrm{in}$ and $V^\mathrm{out}$  respectively, and compare it with the volume of the bounding box, $V^\mathrm{box}$, leading to the relative volume error 
 $\epsilon_V =  \left| V^\mathrm{in} + V^\mathrm{out} - V^\mathrm{box} \right| /  V^\mathrm{box}$.

\begin{figure}[http]
  \begin{subfigure}[b]{.24\textwidth}
    \includegraphics[width=\textwidth]{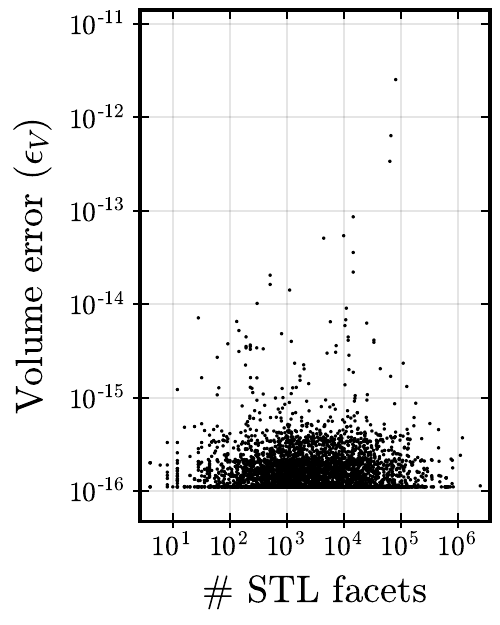}
    \caption{}
    \label{fig:thinki10k-vol-err-nf}
  \end{subfigure}
  \begin{subfigure}[b]{.24\textwidth}
    \includegraphics[width=\textwidth]{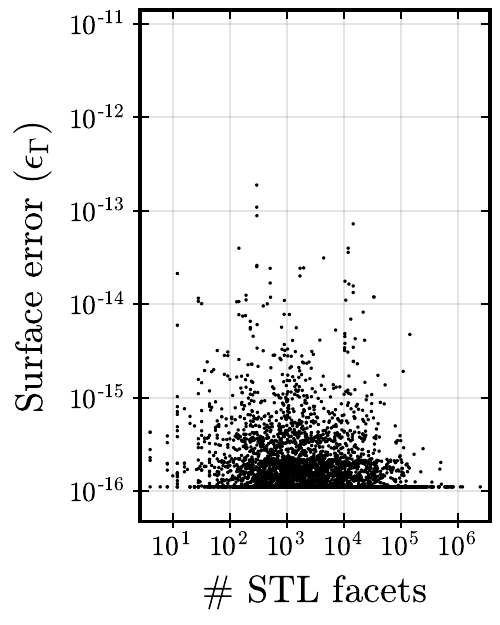}
    \caption{}
    \label{fig:thinki10k-surf-err-nf}
  \end{subfigure}  
  \begin{subfigure}[b]{.24\textwidth}
    \includegraphics[width=\textwidth]{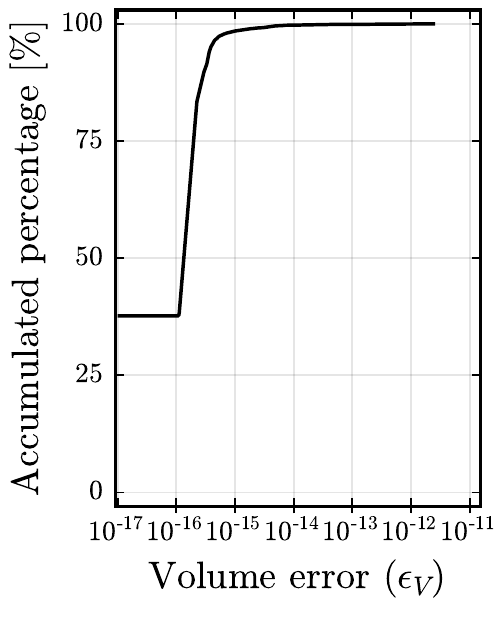}
    \caption{}
    \label{fig:thinki10k-vol-hist}
  \end{subfigure}
  \begin{subfigure}[b]{.24\textwidth}
    \includegraphics[width=\textwidth]{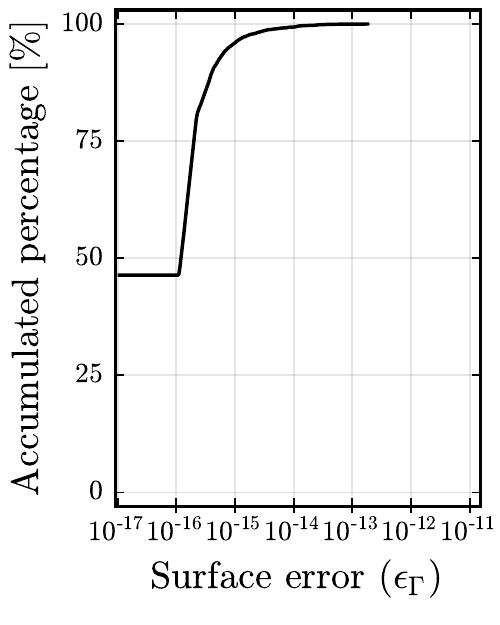}
    \caption{}
    \label{fig:thinki10k-surf-hist}
  \end{subfigure}
  \caption{Volume and surface error distributions: (a) and (b) shows volume and surface errors vs the number of \ac{stl} faces. Each single dot represents a geometry. The cumulative frequency of the volume and surface errors are represented in (c) and (d) respectively.}
  \label{fig:thingi10k-err}
\end{figure}

\fig{fig:thingi10k-err} reports the computed errors $\epsilon_\Gamma$ and $\epsilon_V $ for all processed geometries. Note that the intersection algorithm is able to successfully finish in all cases with relative volume and surface errors below $10^{-11}$ and $10^{-12}$ respectively, which confirms that the method is able to capture the given \ac{stl} models accurately. Note also that the computed errors do not depend on the number of \ac{stl} faces, even for geometries with millions of faces (see  \fig{fig:thinki10k-vol-err-nf} and~\ref{fig:thinki10k-surf-err-nf}). In addition, the volume and surface errors $\epsilon_\Gamma$ and $\epsilon_V $ are below $10^{-15}$ for the virtual majority of cases (see \fig{fig:thinki10k-vol-hist} and \ref{fig:thinki10k-surf-hist}), which demonstrates that the algorithms are able to capture the given \ac{stl} geometries exactly up to the tolerances as expected. Taking into account the 
large number and variety of \ac{stl} models considered, the results of this experiment clearly show that the proposed intersection \alg{alg:global} is able to deal with complex and arbitrary data automatically and provide volume and surface triangulation that capture the input \ac{stl} exactly up to tolerances and close to machine precision.

\subsection{Robustness test}
\label{sec:robustness}

In this second experiment, we study  a sub-set of the models in the Thingi10K database in more detail to assess the robustness of the proposed method with respect to perturbations in the  background mesh. We consider the \ac{stl} geometries displayed in \fig{fig:mat-geo} plus a toy \ac{stl} model of a cube  that will serve as a reference. These \ac{stl} geometries are specifically chosen to cover a large range of shapes and model sizes, while keeping the number of considered cases relatively small in order to make feasible the computation of this example with the computational resources we have at hand.  \tab{tab:test-geos} contains a summary of the main features of the analysed geometries.

\begin{table}[http]
\begin{small}
  \begin{tabular}{rrrcc}
     \toprule
     Model id & Num faces & Num vertices & Box Size & Surface \\
    \midrule
    \texttt{252119} & 49950 & 24979 & (65.06, 37.371, 111.76) & 12439.27 \\
    \texttt{293137} & 292 & 148 & (108.12, 86.625, 107.26) & 29490.72 \\
    \texttt{35269} & 40246 & 20125 & (92.951, 93.426, 33.648) & 8849.629 \\
    \texttt{37266} & 29472 & 14738 & (56.527, 52.541, 53.059) & 13112.22 \\
    \texttt{37881} & 3400 & 1700 & (33.504, 33.688, 18.5) & 3992.616 \\
    \texttt{441708} & 112402 & 56203 & (107.75, 87.802, 107.89) & 29684.95 \\
    \texttt{47076} & 1532 & 768 & (93.095, 93.095, 42.0) & 21864.75 \\
    \texttt{550964} & 6156 & 3072 & (20.0, 20.0, 45.0) & 1715.988 \\
    \texttt{551021} & 348128 & 174066 & (38.365, 25.963, 37.438) & 7282.98 \\
    \texttt{65904} & 157726 & 78869 & (532.5, 552.51, 490.47) & 773637.9 \\
    \texttt{96457} & 1634 & 813 & (195.64, 120.13, 20.549) & 21396.49 \\
    \texttt{cube} & 12 & 8 & (1.0, 1.0, 1.0) & 6.0 \\
    \bottomrule
  \end{tabular}
\end{small}

\vspace{1em} 
 
  \caption{Main features of the test geometries considered displayed in \fig{fig:mat-geo}.}
  \label{tab:test-geos}
\end{table}

The setup of this experiment is as follows. For each \ac{stl} model, we generate different background meshes by perturbing an initial grid, either using translations or rotations.  The initial (unperturbed) mesh for a given \ac{stl} is generated as in previous experiment, but now taking $n^{\rm max}=112$. This value is chosen to stress the algorithm for the reference cube geometry since it leads to faces of the background mesh to be exactly aligned with the faces of the \ac{stl} model. In this scenario, small perturbations of the background mesh lead to volume sub-triangulations with arbitrary small cells, which is a challenging degenerated case. In this regards, we want to analyse how the method behaves, when the perturbation magnitude approaches the machine precision. The first perturbation strategy is to apply a prescribed translation in all directions with magnitude $(\mathbf{p}^{\rm max}_{\rm msh}-\mathbf{p}^{\rm min}_{\rm msh})\Delta x$, where $\Delta x$ is the perturbation coefficient computed as $\Delta_x=10^{-\alpha}$ with $\alpha=1,\ldots,17$. The second perturbation strategy is an imposed rotation composed by three individual rotations of angle $\Delta_\theta$, one over each Cartesian axis, taking the \ac{stl} bounding box barycentre as the origin. 
 
The perturbation angle  is  $\Delta_\theta=10^{-\alpha}$ with $\alpha=1,\ldots,17$. As a result, we consider 34 different background meshes (17 translated + 17 rotated) for each of the \ac{stl} models considered in this example. Finally, we run the intersection  \alg{alg:global} and compute the resulting volume and surface errors  with respect to the non perturbed state $\epsilon_{\Gamma_0}$ and $\epsilon_{V_0} $, which are defined as $\epsilon_{\Gamma_0} = |\Gamma - \Gamma_0| / \Gamma_0$ and $\epsilon_{V_0} = |V - V_0| / V_0$, where  $\Gamma$ and $V$ are the respective surface at each point, and  $\Gamma_0$ and $V_0$ are the respective surface and volume computed with the unperturbed mesh. In contrast to Figure~\ref{fig:thingi10k-err}, here we can take advantage a reference volume.

\begin{figure}[ht!]
  \centering
  \includegraphics[width=0.6\textwidth]{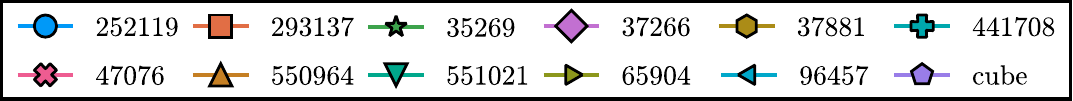}
  
\begin{subfigure}{0.24\textwidth}
\includegraphics[width=\textwidth]{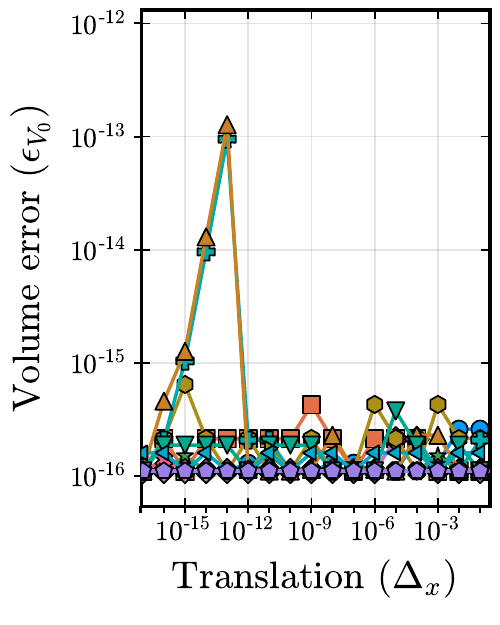}
\caption{}
  \label{fig:robust-trans-vol}
\end{subfigure}
\begin{subfigure}{0.24\textwidth}
\includegraphics[width=\textwidth]{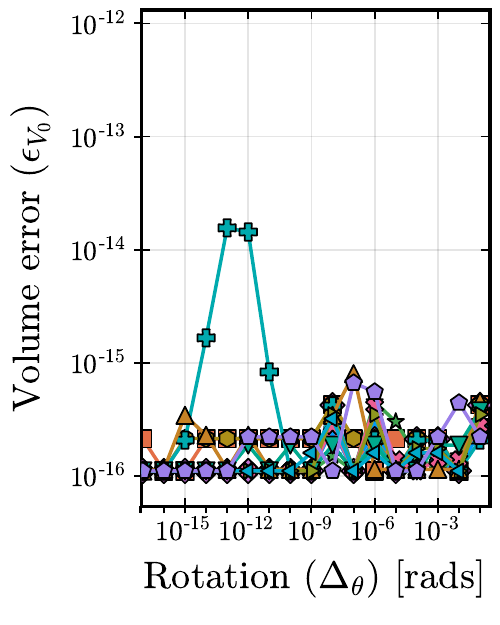}
\caption{}
\end{subfigure}
\begin{subfigure}{0.24\textwidth}
\includegraphics[width=\textwidth]{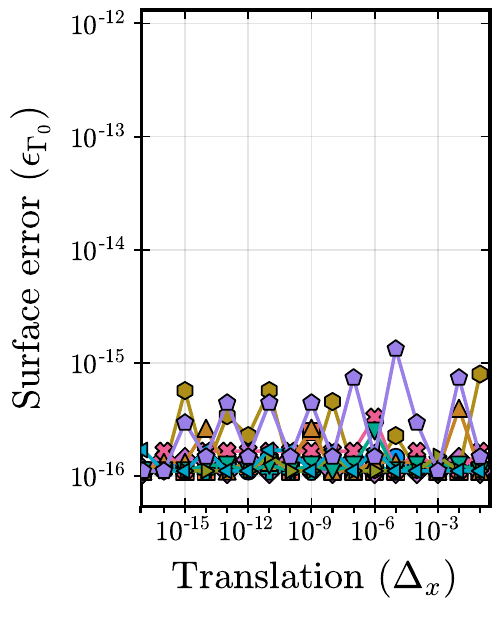}
\caption{}
\end{subfigure}
\begin{subfigure}{0.24\textwidth}
\includegraphics[width=\textwidth]{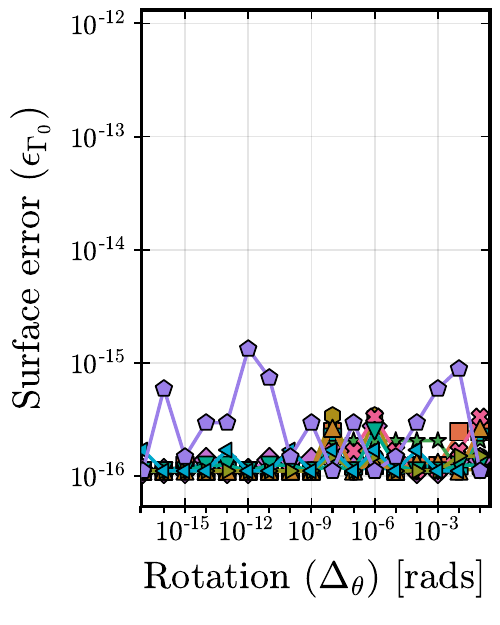}
\caption{}
\end{subfigure}
  \caption{Results of the robustness test: Volume and surface errors $\epsilon_{\Gamma_0}$ and $\epsilon_{V_0} $ in function of the perturbation coefficients $\Delta_x$ and $\Delta_\theta$ for all the geometries of Figure~\ref{fig:mat-geo}.} 
\label{fig:robust}
\end{figure}

As displayed in \fig{fig:robust}, the volume and surface errors $\epsilon_{\Gamma_0}$ and $\epsilon_{V_0} $ are nearly independent on the perturbation coefficients and are below $10^{-15}$ in almost all cases. Some outliers show some influence on the perturbation coefficients but the maximum volume and surface errors are below $10^{-13}$, which is still close to the machine precision and can be attributed to propagation of round-off errors an the value of the tolerance $\epsilon_{hs}$. 
Note that the errors for the cube geometry are always below $10^{-15}$ even though this test has been explicitly designed to render very pathological cases, when the perturbation coefficients tend to zero. At the view of these results, one can conclude that the quality of the computed sub-meshes is nearly independent to the location of the background mesh and, thus, the method is robust to perturbations. 

\subsection{Finite Element convergence test}
\label{sec:fe-test}

 In this last experiment, we explore the capacity of the proposed intersection algorithm to be coupled with unfitted \ac{fe} methods in order to simulate complex geometries described by \ac{stl} models without generating conforming unstructured grids. 
  The main goal of this experiment is to check that the intersection algorithm does not introduce any spurious numerical artefacts that destroy the optimal convergence of the \ac{fe} solver. We will also leverage this convergence test to evaluate the performance of the intersection method by studying the scaling of CPU times with respect to the number of cells in the background mesh.
 
For the \ac{fe} computation, we consider a Poisson equation with pure Dirichlet boundary conditions as the model problem. 
A numerical approximation $u_h\approx u$ is computed with the \ac{agfem} method described in~\cite{Badia2018c} for exactly the same model problem. In particular, the interpolation spaces are defined with continuous tri-linear Lagrangian shape functions.  As an example, see in  \fig{fig:sin-functions}, \ac{fe} approximations computed on a sub-set of the studied \ac{stl} models. 
 
\begin{figure}[http]
  \includefig{0.4\textwidth}{legend_poisson}

  \begin{subfigure}[b]{0.35\textwidth}
    \includegraphics[width=\textwidth]{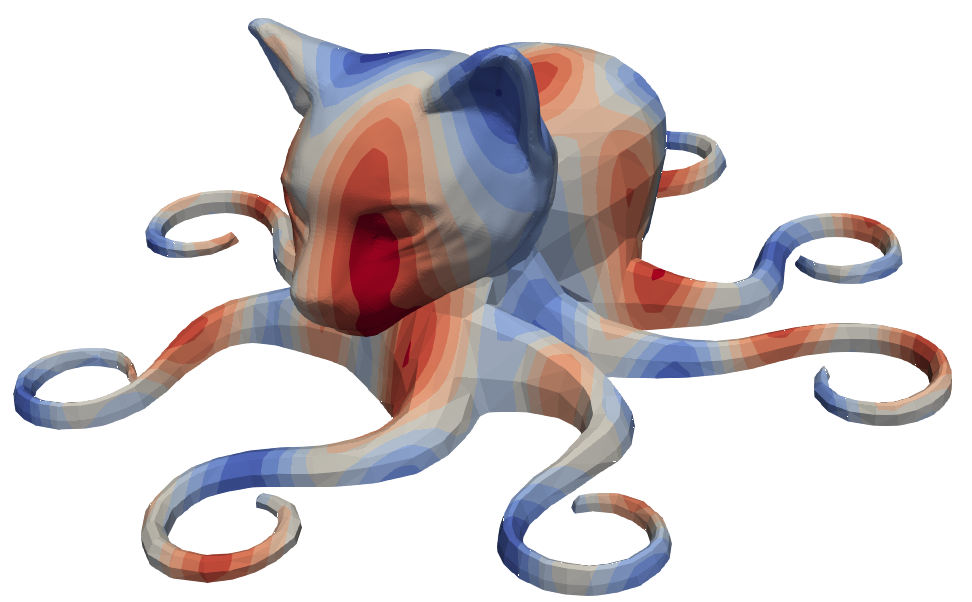}
    \caption{35269}
  \end{subfigure}
  \begin{subfigure}[b]{0.2\textwidth}
    \includegraphics[width=\textwidth]{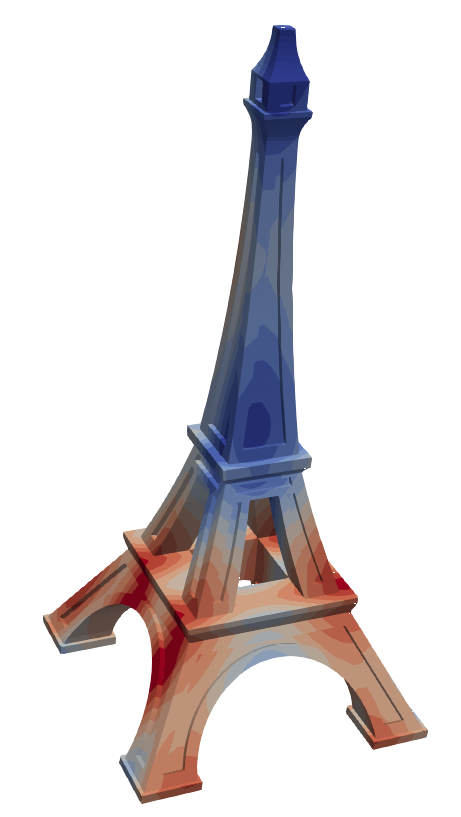}
    \caption{55094}
  \end{subfigure}
  \begin{subfigure}[b]{0.35\textwidth}
    \includegraphics[width=\textwidth]{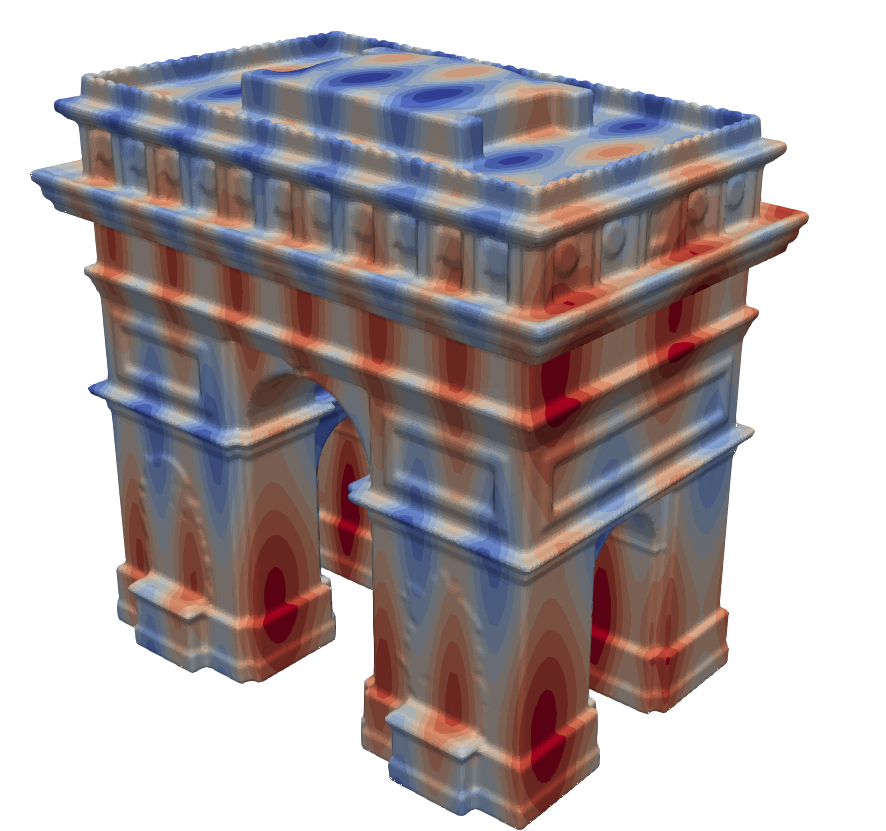}
    \caption{551021}
  \end{subfigure}
  \caption{\ac{fe} approximation computed with \ac{agfem} on top of three of the \ac{stl} models analysed in the experiments. Here, the underlying Poisson equation is defined using the manufactured solution $u(x,y,z) = sin(a \tfrac{ 2\pi }{ T }x )+ sin(b \tfrac{ 2\pi }{ T } y) + sin(c \tfrac{ 2\pi }{ T } z)$ with $T = 10 h$, $(a,b,c) = (1,\tfrac{1}{2},\tfrac{1}{4})$, 
  $h = \frac{1.4}{n^{\mathrm{max}}} \max( \mathbf{p}^{\rm max}_{\rm stl}-\mathbf{p}^{\rm min}_{\rm stl} ) $  
  and $n^{\mathrm{max}}=100$.}
  \label{fig:sin-functions}
\end{figure}

For the convergence test, the forcing term and Dirichlet boundary condition are defined such that the manufactured function $u(x,y,z) = x^2+y^2-z^2$ is the exact solution of the problem. Since this function is smooth and does not belong to the interpolation space, we expect that the $H^1$ and $L^2$ norms of the discretisation error $e_h\doteq u-u_h$, namely
\begin{equation}
\|e_h\|^2_{L^2(\Omega)} \doteq \int_{\Omega} e_h^2 {\rm\ d}\Omega \quad \text{ and }\quad
\|e_h\|^2_{H^1(\Omega)} \doteq \int_{\Omega} e_h^2 + \nabla e_h\cdot\nabla e_h {\rm\ d}\Omega, 
\end{equation}
converge with the optimal convergence rate. Our goal is to compute these error norms for different mesh sizes and confirm that they converge with the optimal slopes. 

We build a family of background meshes for each \ac{stl} model in \fig{fig:mat-geo}. Each mesh is generated by using a different value of $n^{\rm max}$, leading to several refinement levels. In particular, we use $n^{\rm max}=n^{\rm max}_0 2^\beta$ with $n^{\rm max}_0=14$ and  $\beta=0,\ldots,5$. The finest meshes generated in this way ($\beta=5$) have $448$ cells in the largest axis. In order to be able to solve the underlying system of linear algebraic equations for such problem sizes, we consider a conjugate gradient solver preconditioned with the \ac{amg} method in the Preconditioners.jl package version 0.3} \cite{Preconditioners-jl}. We declare convergence of the conjugate gradient solver, when the relative energy norm is below $10^{-10}$. 

\begin{figure}[http]
  \centering
  \includegraphics[width=0.833\textwidth]{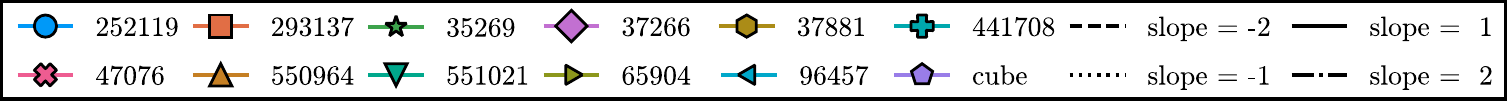}
  
  \begin{subfigure}[b]{.24\textwidth}
    \includegraphics[width=\textwidth]{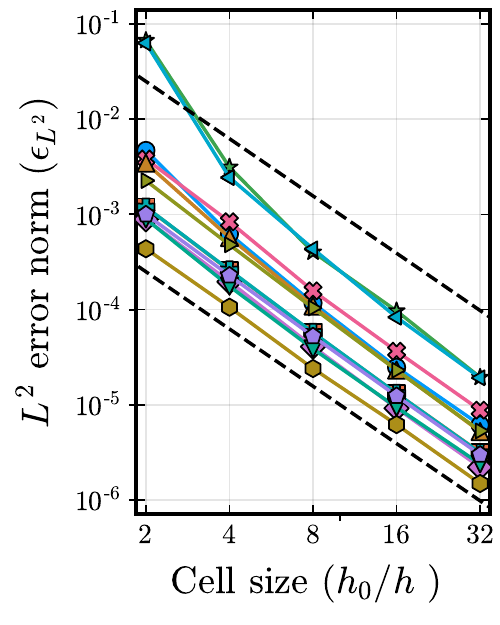}
    \caption{}
  \end{subfigure}
  \begin{subfigure}[b]{.24\textwidth}
    \includegraphics[width=\textwidth]{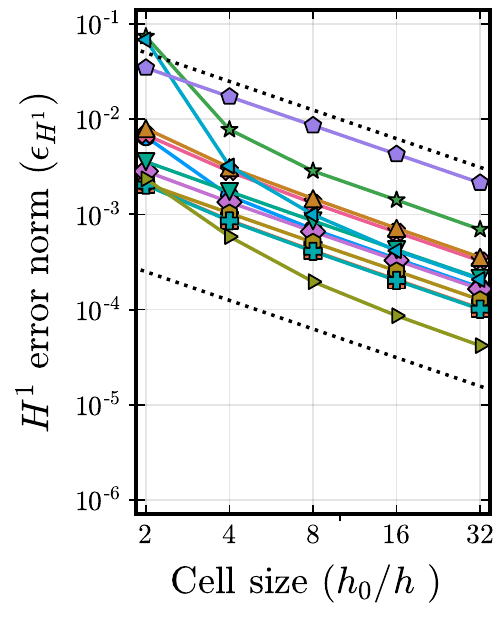}
    \caption{}
  \end{subfigure}
    \begin{subfigure}[b]{0.24\textwidth}
    \includegraphics[width=1\textwidth]{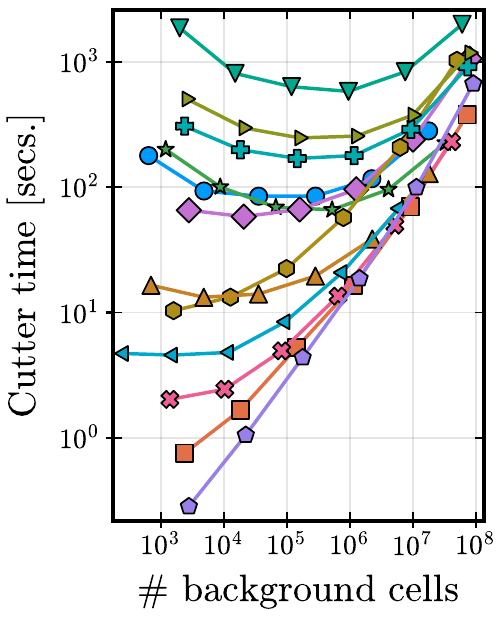}
    \caption{}
    \label{fig:agfem-conv-3}
  \end{subfigure}
  \begin{subfigure}[b]{0.24\textwidth}
    \includegraphics[width=1\textwidth]{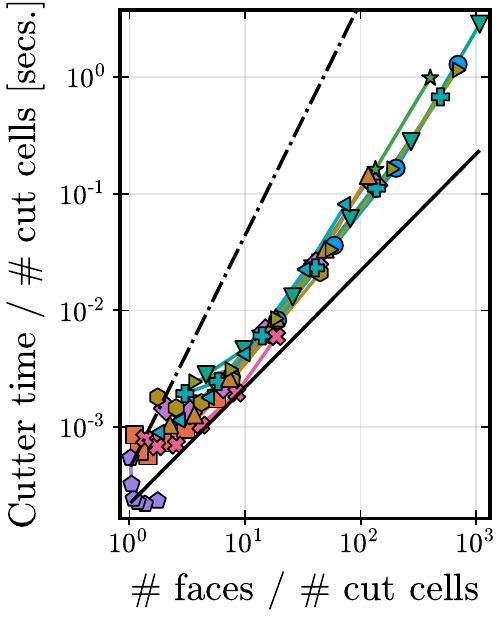}
    \caption{}
    \label{fig:agfem-conv-4}
  \end{subfigure}
  \caption{ Results of the \ac{fe} convergence test: (a), (b), and (c) show the  $L^2$, and $H^1$ error norms, and CPU times in \alg{alg:global} vs relative cell size for all \ac{stl} models in \fig{fig:mat-geo}. (c) shows the scaling of CPU times vs number of \ac{stl} faces per cut cell in the background mesh.}
  \label{fig:agfem-conv}
\end{figure}

The results of the convergence test are displayed in \fig{fig:agfem-conv}. The $L^2$ and $H^1$ error norms converge with the expected slopes for all the considered \ac{stl} geometries, which confirms that the intersection algorithm is not affecting the quality of the \ac{fe} solver. On the other hand, we measure the CPU time elapsed in the computation of the intersection \alg{alg:global}.  \fig{fig:agfem-conv-3} shows the scaling of the CPU time with respect to the number of cells in the background mesh. The scaling tends to be linear as the mesh is refined in all cases. The linear regime is reached at different speeds depending on the considered \ac{stl} models. For the cube, which is the simplest geometry studied here, the linear regime is reached before the other ones, whereas the model with more \ac{stl} faces (the \emph{Arc de Triomphe} geometry with id 551021) is the latest one to achieve the linear regime. Note that the CPU times converge to similar values for all geometries, when the mesh is refined. This is because, in the limit, the cut algorithm only needs to intersect each cut background cell with a single plane independently of the number of faces in the \ac{stl}. This suggest that the number of \ac{stl} planes per cut cell in the background mesh is closely related with the performance of the method. To analyse the interplay between these two quantities, \fig{fig:agfem-conv-4} displays the scaling of CPU time with respect to the number of \ac{stl} faces both averaged by the number of cut cell in the background mesh. Due to the nature of the algorithm, which involves searches between the \ac{stl} and the background mesh, one can expect a superlinear scaling. This is indeed what is observed in \fig{fig:agfem-conv-4}, but, in any case, the scaling is clearly not quadratic since the searches are efficiently computed using a tree partition. In particular,  this allowed us to compute sub-triangulations of complex \ac{stl} geometries with hundreds of thousands of \ac{stl} faces in this example and even millions of \ac{stl} faces in previous examples in \sect{sec:thingi10k}.

\begin{figure}[http]
  \begin{subfigure}[b]{0.49\textwidth}
    \raggedright
    \includegraphics[width=0.95\textwidth]{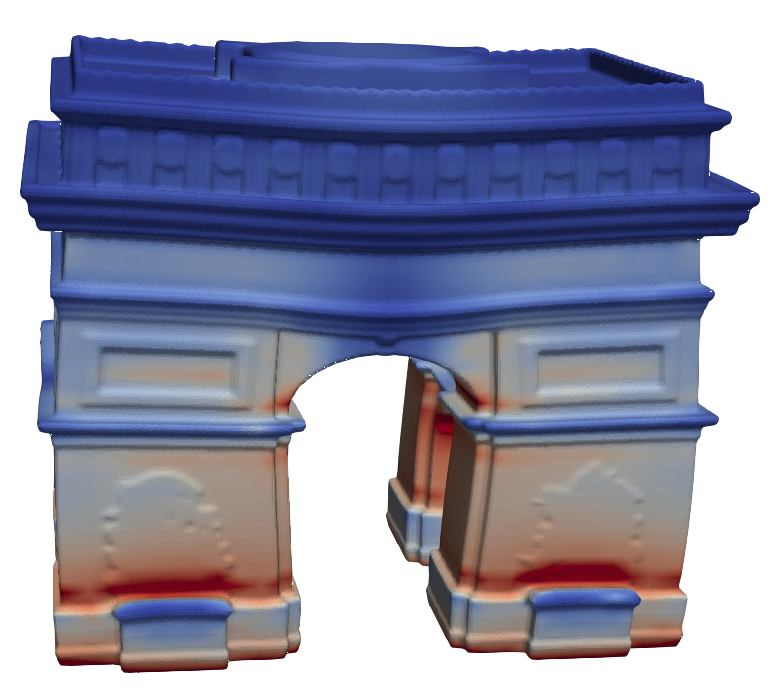}
    \caption{}
  \end{subfigure}
  \begin{subfigure}[b]{0.49\textwidth}
    \includegraphics[width=\textwidth]{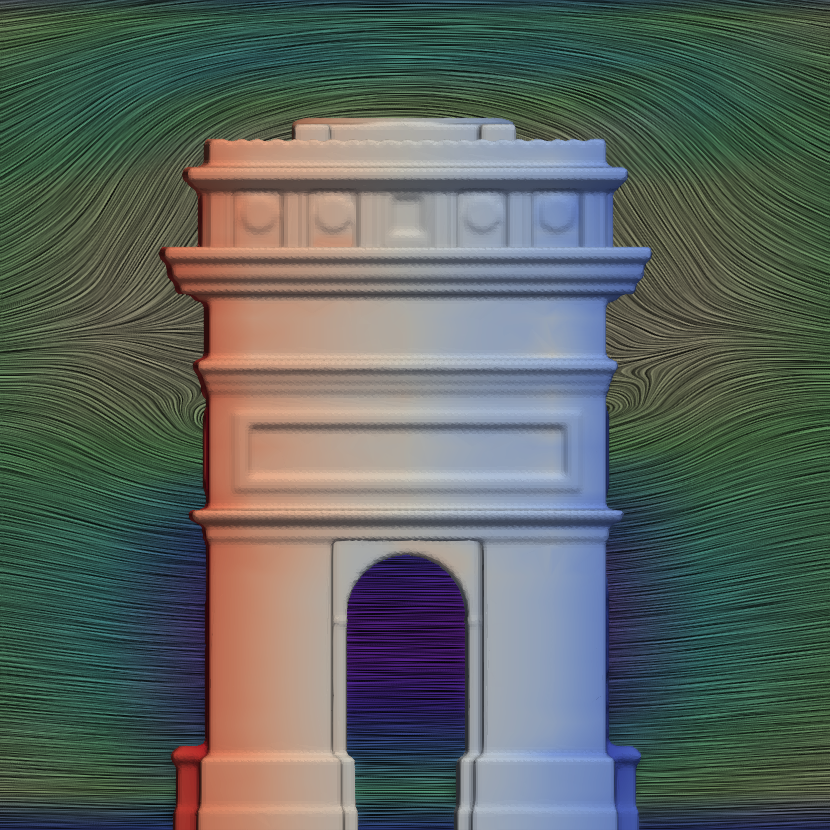}
    \caption{}
  \end{subfigure}
  \caption{We show the {\ac{agfem} approximations of two physical problems on both sides of the \emph{Arc de Triomph} \ac{stl}. In (a), we show the deformed configuration for linear elasticity and the colour map for the stress field. A vertical body force $(0,0,-1)$ is applied to the volume, representing its own weight. The elastic modulus is set to $10^{-5}$, the Poisson ration is 0.3 and the deformation is magnified 500 times.
  We use a Cartesian mesh with  $50\times50\times50$ cells.
  In (b), we show a line integral convolution of the velocity field for an inviscid and incompressible flow around the geometry and the pressure colour map on the surface. The inlet velocity is set as $v_\mathrm{in}(x,y,z) = v_\mathrm{max}\left(0,v_1(x) v_1(z),0\right)$, where $v_\mathrm{max} = 0.2$ and $v_1(x) = 4x-4x^2 $. The wall velocity is zero.  We use a Cartesian mesh with $20\times60\times20$ cells.} }
  \label{fig:phys-examples}
\end{figure}

{The proposed geometrical treatment can readily be applied to other unfitted \acp{fe} methods and \acp{pde}. In Figure~\ref{fig:phys-examples}, we solve a linear elasticity problem in the \emph{Arc de Triomphe} geometry and an incompressible flow problem surrounding it. For linear elasticity, we used the formulation in \cite{Neiva2021}. The unfitted method for incompressible flows can be found in \cite{Badia2018a}. One can also observe in these two examples that we can readily use the meshes on both sides of the boundary representation.}

\section{Conclusions and future work}\label{sec:conclusions}
 
In this work, we have designed a fully automatic simulation pipeline for the numerical approximation of \acp{pde} on general domains described by a boundary mesh. The algorithm makes use of a structured background mesh and an unfitted \ac{fe}
formulation on this mesh. The main complication of these methods is the numerical computation of integrals in the interior of the domain for background cells cutting the domain boundary. Boundary meshes for complex geometries can involve a huge number of faces intersecting background cells and the geometries are not convex in general.

We have designed a general clipping algorithm for cut cells that can deal with general surface meshes. They are based on convex decomposition algorithms, robust clipping of convex polyhedra, a graph-based representation of polyhedra, discrete level-set representation of planes and some merging techniques to reduce rounding error effects. The result of this algorithm is a refinement of the boundary mesh that can readily be used to integrate boundary terms and a two-level \emph{integration} mesh. The two-level mesh combines the background mesh and a cell-wise partition of cut cell interiors into convex polyhedra and can straightforwardly be used to integrate the bulk terms in unfitted \ac{fe} schemes.

The implementation of the algorithms are distributed as open source software and can be found in \cite{Martorell_STLCutters_2021}. The algorithm implementation has been applied with success on all 3D analysis-suitable  meshes in the Thingi10K database~\cite{Zhou2016} (almost 5,000 meshes), showing the sound robustness of the approach. The reported integration errors are close to machine precision, which prove its accuracy. These integration meshes have been successfully combined with one unfitted formulation, the \ac{agfem} \cite{Badia2018c}, to discretise \acp{pde} on these geometries, and convergence error plots are provided. Finally, the computational complexity and cost of the geometrical algorithm is reported and compared against the \ac{fe} solver step.

Future work involves the extension of this approach to other background meshes, specially octree meshes, even though this extension is quite straightforward; the algorithms are cell-wise defined and can readily be applied to locally refined structured meshes. It is of practical relevance to extend the current (open source) implementation to distributed-memory computers. Since the algorithms mainly involve cell-wise computations, they are embarrassingly parallel, a great benefit compared to unstructured mesh generation algorithms that require global consistency and thus are very hard to parallelise. The extension to 4D (under homotopy assumptions or allowing topology changes in time) is of special relevance, since it would allow one to solve complex time-dependent problems (e.g., fluid-structure interaction or multi-fluid models) that involve moving interfaces, one of the main challenges in the field. Most of the ingredients in the current algorithms are dimension-agnostic and the polyhedron representation in terms of oriented graphs seems to be general enough (the formulation does not rely on planar graphs, which would prevent a 4D extension, since 4D polytopes cannot be represented as planar graphs in general). Another topic of interest is the extension of this approach to higher order boundary representations, e.g., connecting the algorithm with the B-REP representation to attain higher levels of accuracy via nonlinear intersection algorithms.

\section*{Acknowledgments}

\newcommand{\thethanks}{This research was partially funded by the Australian Government through the Australian Research Council (project number DP210103092), the European Commission under the FET-HPC ExaQUte project (Grant agreement ID: 800898) within the Horizon 2020 Framework Programme and the project RTI2018-096898-B-I00 from the ``FEDER/Ministerio de Ciencia e Innovación – Agencia Estatal de Investigación''. F. Verdugo acknowledges support from the Spanish Ministry of Economy and Competitiveness through the ``Severo Ochoa Programme for Centers of Excellence in R\&D (CEX2018-000797-S)". P.A. Martorell aknowledges the support recieved from Universitat Politècnica de Catalunya and Santander Bank through an FPI fellowship (FPI-UPC 2019). This work was also supported by computational resources provided by the Australian Government through NCI under the National Computational Merit Allocation Scheme.}
 
\thethanks

\setlength{\bibsep}{0.0ex plus 0.00ex}
\bibliographystyle{myabbrvnat}
\bibliography{refs} 
  
\end{document}